\documentclass[11pt]{amsart}
\usepackage{amsfonts,amssymb,amsthm,amsmath,amsxtra,amscd,verbatim,eucal}
\usepackage[all]{xy}
\usepackage[dvips]{graphics}
\usepackage{graphicx}

\theoremstyle{plain}
\newtheorem{thm}{Theorem}[section]

\newtheorem*{thmA}{Theorem~A}
\newtheorem*{thmD}{Theorem~D}
\newtheorem*{corB}{Corollary~B}
\newtheorem*{corC}{Corollary~C}

\newtheorem{lem}[thm]{Lemma}
\newtheorem{prop}[thm]{Proposition}
\newtheorem{cor}[thm]{Corollary}

\theoremstyle{definition}

\theoremstyle{remark}
\newtheorem{eg}[thm]{Example}

\newtheorem{rmk}[thm]{Remark}

\numberwithin{equation}{section}

\def\Z{{\mathbb Z}}

\def\C{{\mathbb C}}
\def\R{{\mathbb R}}
\def\Q{{\mathbb Q}}
\def\P{{\mathbb P}}
\def\bB{{\mathbb B}}

\def\CC{\mathcal{C}}

\def\K{\mathcal{K}}
\def\O{\mathcal{O}}

\def\frD{\mathfrak{D}}

\def\a{\alpha}
\def\b{\beta}

\def\d{\delta}
\def\f{\phi}

\def\l{\lambda}

\def\s{\sigma}

\def\D{\Delta}

\def\Si{\Sigma}

\def\.{\cdot}
\def\~{\widetilde}
\def\^{\widehat}

\def\o{\circ}

\renewcommand{\and}{ \quad \text{and} \quad }

\DeclareMathOperator{\graff} {Graff}

\DeclareMathOperator{\Tr}  {Tr}

\DeclareMathOperator{\Spec} {Spec}

\DeclareMathOperator{\Pic} {Pic}

\DeclareMathOperator{\End} {End}

\DeclareMathOperator{\Val} {Val}

\DeclareMathOperator{\Vol} {Vol}

\DeclareMathOperator{\id} {id}

\DeclareMathOperator{\Hom} {Hom}

\DeclareMathOperator{\gr} {Gr}
\DeclareMathOperator{\pl}{Pl}

\def\rO{\mathrm{O}}
\def\cD{\mathcal{D}}

\def\bG{{\mathbb G}}
\def\dto{\dashrightarrow}
\DeclareMathOperator{\Hproj}{{\mathcal H}^* }
\DeclareMathOperator{\hproj}{{\mathcal H}}
\DeclareMathOperator{\Hinj}{\underrightarrow{\mathrm H}^*}
\DeclareMathOperator{\hinj}{\underrightarrow{\mathrm H}}
\DeclareMathOperator{\PL}{{\mathrm PL}}
\DeclareMathOperator{\PP}{{\mathrm PP}}
\DeclareMathOperator{\ch}{{\mathrm ch}}

\DeclareMathOperator\tch{{\mathrm coh}}
\def\1{{\bf 1}}
\def\aff{\mathrm{aff}}
\def\Atrans{\check A}
\def\A{A}
\def\B{B}
\DeclareMathOperator{\dVol} {dVol}

\title[Degree growth of monomial maps]{Degree growth of monomial maps \\ and McMullen's polytope algebra}

\date{\today}
\thanks{First author is supported by the ANR-project Berko. Second author partially supported by the Swedish Research Council and the NSF}

\author[C. Favre]{Charles Favre}
\address{Centre de Math\'ematiques Laurent Schwartz\\
\'Ecole Polytechnique\\
91128 Palaiseau Cedex \\France}
\email{favre@math.polytechnique.fr}

\author[E. Wulcan]{Elizabeth Wulcan}
\address{Department of Mathematics\\Chalmers University of Technology and the University of G\"oteborg\\S-412 96 G\"oteborg\\SWEDEN}
\email{wulcan@chalmers.se}

\thanks{}

\begin{document}

\begin{abstract}
We compute all dynamical degrees of monomial maps by interpreting them as 
mixed volumes of polytopes. By exploiting further the 
isomorphism between the polytope algebra of P.~McMullen and the universal 
cohomology of complete toric varieties, we construct invariant positive 
cohomology classes when the dynamical degrees have no resonance.
\end{abstract}

\maketitle

\section{Introduction}
Some of the most basic information associated to a rational dominant map $f: \P^d \dashrightarrow \P^d$ is provided by its degrees $\deg_k(f) := \deg f^{-1} (L_k)$, where $L_k$ is a generic
linear subspace of $\P^d$ of codimension $k$, see~p.\pageref{def-pullback} below for a formal definition.  From a dynamical point
of view, it is important to understand the behaviour of the sequence
$\deg_k(f^n)$ as $n\to\infty$. 
It is not difficult to see that $\deg_k(f^{m+n }) \le\deg_k(f^m)
\,\deg_k(f^n)$, and thus following Russakovskii-Shiffman ~\cite{RS} we
can define the \emph{$k$-th 
dynamical degree of $f$} as $\l_k(f) :=\lim_n \deg_k(f^n)^{1/n}$. 
Basic properties of dynamical degrees can be found in~\cite{RS,DS}. 
Our main objective is to describe the sequence of degrees $\deg_k(f^n)$ in the special case of monomial maps $f$, but for arbitrary $k$. 

Controlling the degrees of iterates of a rational map is a quite delicate problem. 
Up to now, most investigations have been focused on the case
$d=2$ and $k=1$, see~\cite{DF,FJ} and the references therein. 
There are also various interesting families of examples for $k=1$ in
arbitrary dimensions in e.g.~\cite{AABM,AMV,BK1,BK3,BHM,N}. In particular, the case of monomial maps and $k=1$ 
is treated  in~\cite{BK2,Fa,HP,JW,L}. 
On the other hand, there are only few references in the literature concerning  the 
case $2\le k \le d-2$, see ~\cite{Og,DN}.  An essential problem arises from the difficulty  
to  explicitly compute $\deg_k(f)$ even in concrete examples.  
This can be 
overcome in the case of monomial maps, 
 since tools from convex geometry allow 
one to compute these numbers in terms of (mixed) volumes of polytopes. This technique has already been used to compute all degrees of the standard Cremona transformation in arbitrary dimensions by  Gonzalez-Springer and Pan~\cite{GSP}.

\smallskip

Monomial maps on $\P^d$ correspond to integer valued $d\times d$ matrices, $M(d,\Z)$. Given $\A\in  M(d,\Z)$ we write $\f_\A$ for the corresponding monomial map $\f_\A(x_1, \ldots , x_d) = (\prod_j x_j^{a_{j1}}, \ldots,  \prod_j
x_j^{a_{jd}})$ with  $(x_1, ...,x_d)\in (\C^*)^d$. 
This mapping is holomorphic on the torus $(\C^*)^d$ and extends
as a rational map to the standard equivariant compactification $\P^d
\supset (\C^*)^d$. Moreover $\f_\A$ is dominant precisely if $\det (\A)\neq 0$. 
Observe that $\f_\A^n = \f_{\A^n}$ for all $n$.
If $a_n$ and $b_n$ are sequences of positive real numbers, we write
$a_n\asymp b_n$ if 
$C^{-1}\le a_n/b_n \le C$ for some $C>1$ and all $n$. 
\begin{thmA}\label{thm:main1}
Let $\A\in M(d,\Z)$
 and let $\f_\A: \P^d \dto \P^d$ be the corresponding rational map.  Assume that $\det (\A) \neq 0$ (so that $\f_\A$ is dominant). 
Then, for $0\le k \le d$, 
\begin{equation}
\deg_k(\f_\A^n) \asymp \| \wedge^k \!\!\A^n\|,
\end{equation}
where $\wedge^k \!\A : \wedge^k \R^d \to \wedge^k \R^d $ is the natural linear map induced by $\A$ and 
$\|\cdot\|$ is any norm on $\End(\wedge^k \R^d)$. 
\end{thmA}
\begin{corB}\label{cor:dyndegree}
Let $\f_\A$ and $\A$ be as in Theorem A. 
Order the eigenvalues of $\A$ in decreasing order, $|\rho_1| \ge |\rho_2|\ge \ldots  \ge |\rho_d|$.  
Then the  $k$-th dynamical degree of the monomial map $\phi_\A$ is equal to
$\prod_1^k |\rho_j|$. 
\end{corB}

Recall that the topological entropy of a rational map $f : X \dto X$ on a projective variety is defined as the asymptotic rate of growth of $(n,\varepsilon)$-separated sets outside the indeterminacy set of iterates of $\f$, 
see~\cite{DS} for details.
On the one hand, the topological entropy of a monomial map is greater than its restriction 
to the compact real torus
$\{ |x_i| =1\}\subset (\C^*)^d$ which is equal to $\log (\prod_1^d \max \{ 1, |\rho_i|\})$, see~\cite[Sect.~5]{HP}. 
On the other hand, it is a general result due to Gromov~\cite{G}, and Dinh-Sibony~\cite{DS} 
that $\max_k \log \l_k$ is an upper bound for the topological entropy. By Corollary B, $\log (\prod_1^d \max \{ 1, |\rho_i|\})=\max_k \log \l_k$.
Thus we have 
\begin{corC}\label{cor:entropy}
Let $X$ be a projective smooth toric variety, let $\A$ be as in Theorem A, and let $\f_\A:X\dto X$ be the induced rational map. Then the topological entropy of $\f_\A$ is 
equal to $\max \log \l_k$.
\end{corC}
We note that Theorem A and its two corollaries have been obtained independently by
Jan-Li Lin, \cite{L2}, by different but related methods. His approach relies on the notion of Minkowski weights. 

\medskip

By Khovanskii-Teissier's inequalities,  the sequence $k\mapsto \log \deg_k(f)$ is concave so that we always have
$\l_k^2(f) \ge \l_{k-1}(f) \, \l_{k+1}(f)$ for any $1\le k \le d-1$.
Our next result gives a more precise control of the degrees when the
asymptotic degrees are strictly concave. It can be seen as an analogue of~\cite[Main Theorem]{BFJ} in the case of monomial maps but in arbitrary dimensions.
\begin{thmD}\label{thm:main2}
Let $\A\in M(d,\Z)$ and let $\f_\A: \P^d\dto
\P^d$ be the associated rational monomial map. Write
$\l_k=\l_k(\f_\A)$. Assume that $\det(\A) \neq 0$ and that for some $1\le k \le d-1$ the dynamical degrees satisfy
\begin{equation}\label{strikt}
\l_k^2> \l_{k-1} \, \l_{k+1}~.
\end{equation}
Then there exists a constant $C>0$ and an integer $r\ge0$  such that, 
for this $k$, 
\begin{equation}\label{lor}
\deg_k(\f_\A^n)  = C \l_k^n + \O \left(n^r\, \left(\frac{\l_{k-1} \, \l_{k+1}}{\l_k}\right)^{n}\right).
\end{equation}
\end{thmD}
Theorems A and D (and thus Corollary B) hold true for $\P^d$ replaced
by a projective smooth variety, cf. Remark ~\ref{general}.

For $k=1$ Theorem~A is proven in~\cite{HP}, and 
Theorem D is due to Lin \cite[Thms~6.6-7]{L}. In fact, there
 are finer estimates for the growth of $\f_A$. For
example, Bedford-Kim~\cite{BK2} gave a description of when $\deg_1(\phi_A^n)$
satisfies a linear recurrence; in particular, it happens if 
$|\rho_1|>|\rho_2|$. We do not know if the assumption in Theorem D is
sufficient for $\deg_k(\phi_A^n)$ to satisfy a linear recursion. 
This problem is related to the construction of a toric model $X(\D)$
dominating $\P^d$ such that the induced action $\f_A^*: H^k(X(\D), \R)\to
H^k(X(\D), \R)$ of the monomial map $\f_\A$ commutes with iteration;
or, in the terminology of Fornaess-Sibony, a model in which the map
induced by $\f_\A$ is \emph{stable}. This very delicate problem
is treated in detail in various papers by Bedford-Kim~\cite{BK2},
the first author~\cite{Fa} ,
Hasselblatt-Propp~\cite{HP}, Jonsson and the second
author~\cite{JW}, and Lin~\cite{L}  in the case $k=1$. 
We do not address this problem here.

\smallskip

Let us briefly explain the idea of the proofs of Theorems~A and~D. 
The degree $\deg_k(\f_A)$ can be naturally interpreted as an intersection number
$\deg_k(\f_A)=\f_\A^* \O(1)^k\cdot \O(1)^{d-k}$. 
Recall that any polytope $P$ with integral vertices determines a toric variety $X(\D_P)$ and a
line bundle $L_P$ over $X(\D_P)$. Conversely, any line bundle on a toric variety that
is generated by its global sections determines a polytope. 
Using this correspondence, one can compute intersection products of line bundles in a toric variety  in terms of mixed volumes of polytopes, see~\cite[p.79]{Od}.
More precisely,  given any two polytopes $P, Q$ giving rise to two line bundles $L_P,L_Q$ on the same toric variety, then  the intersection product $L_P^k\cdot L_Q^{d-k}$ is given as the
\emph{mixed volume} $d! \Vol (P[k], Q[d-k])$,  i.e.,  (a constant times) the coefficient of
$t^k$ in the polynomial $\Vol(tP + Q)$. 
Since $\O(1)$ over $\P^d$ corresponds to the standard simplex  $\Sigma_d\subset\Q^d$, and the action of $\phi_A^*$ on $L_P$ corresponds to the linear action of $A$ on $P$ it turns out that  
the proofs of Theorems ~A and ~D amounts
to controlling the growth of mixed volumes under the action of the
linear map $\A$:
\begin{equation}\label{volym}
\deg_k(\f_\A^n)  = d! \Vol \left ( \A^n(\Si_d) [k], \Si_d[d-k] \right).
\end{equation}
The computation of mixed volumes is in general quite
difficult. However, since
we are interested in the asymptotic behaviour of $\deg_k$ we may
replace $\Sigma_d$ by a ball, which allows us to apply the
Cauchy-Crofton formula, see Section~\ref{sec:cvx1}. When $A$ is
diagonalizable, one can avoid this geometric-integral formula and work
in a basis of diagonalization of $A$ to estimate directly the growth
of $\deg_k(\f_\A^n)$, see Section~\ref{sec:diag}.

\smallskip

For dynamical applications it is often crucial to construct invariant cohomology classes with
nice positivity properties.
%such as nef classes. 
In Section ~\ref{sec:inv}, we explain how to construct such invariant classes for monomial maps satisfying the assumptions of Theorem~D. Our Corollary~\ref{klasser} is an analog of \cite[Corollary~3.6]{BFJ} in the toric setting and in arbitrary dimensions. 
However these classes do not live in the cohomology groups of a particular toric variety,  but rather on the inductive limit of all cohomology groups
$\hinj^k := \varinjlim H^{2k}(X(\D), \R)$ over \emph{all} toric
models. This idea has already been 
fruitfully used in dynamics in~\cite{C,BFJ}, and we propose a general framework to study
 the action of monomial maps on these spaces.

To this end, we rely on a beautiful interpretation of $\Hinj$ in
terms of convex geometry due to Fulton-Sturmfels 
\cite{FS}, see also~\cite{B}.  
Namely, the classes in $\Hinj$ are
in one-to-one correspondence with the classes in 
P. McMullen's \emph{polytope algebra}. 
The polytope algebra $\Pi$ 
is the $\R$-algebra generated by
classes $[P]$ of polytopes with vertices in $\Q^d$, with relations $[P+v]= [P]$ for
$v\in\Q^d$ and $[P\cup Q] + [P\cap Q] = [P] + [Q]$ whenever $P\cup Q$
is convex. It is endowed with multiplication $[P]\cdot[Q]:=[P+Q]$, where $P+Q$
denotes the Minkowski sum. 
To any polytope $P$, we can attach the Chern character of its associated line bundle
$L_P$, and this  defines a linear map $\ch:
\Pi\to\Hinj$, which is, in fact, an  isomorphism of algebras,
~\cite{B,FS}.  It holds that $\f_A^*\ch [P]=\ch[A(P)]$. 

The isomorphism $\ch:\Pi\to\Hinj$ extends by duality to an isomorphism 
between the space of linear forms on $\Pi$ and $\Hproj$. We will  
call elements in the former space  \emph{currents}. The invariant
cohomology classes alluded to above correspond to a 
 very special type of currents obtained by taking the volume of the
 projection of the polytope on suitable linear subspaces.

Although strictly not needed for the proofs, the  formalism of currents on the polytope algebra is mainly motivated by
our endeavour of constructing invariant classes for monomial maps. 
However, we note that the space of currents contains classical objects from 
 convex geometry, such as valuations in the sense of~\cite{MS}. We
 think that it would be interesting to further explore this space; e.g. investigate positivity properties of currents and 
define (under reasonable geometric conditions) the intersection product of currents.

\smallskip

The paper is organized as follows. 
Sections~\ref{sec:basics} and ~\ref{sec:polytope} contain basics on
toric varieties and the polytope algebra, respectively. 
In Section \ref{sec:dynamical} we discuss dynamical degrees
on toric varieties and, in particular, we derive \eqref{volym}. The
proof(s) of Theorem~A (and Corollary~B) occupies
Section~\ref{sec:pfA}, whereas Theorem~D is proved in Section~\ref{sec:pfD}, and invariant classes of monomial maps are constructed in Section~\ref{sec:inv}.

\smallskip

\noindent
\emph{Acknowledgement:} 
We thank the referee for a careful reading and for many helpful suggestions, and J.-L. Lin for many interesting discussions.

%%%%%%%%%%

\section{Toric varieties}\label{sec:basics}

%%%%
A toric variety $X$ over $\C$ is a normal irreducible algebraic variety endowed with an action 
of the multiplicative torus $\bG_m^d := (\C^*)^d$ which admits an open and dense orbit. 
This section contains the necessary material from toric geometry that will be 
needed for the proof of our results. Our basic references are~\cite{F,Od}.

\subsection{Fans and toric varieties}
Let $N\simeq \Z^d$ be a lattice, i.e.,  a free abelian group, of rank
$d$, denote by $M = \Hom(N,\Z)$ its dual lattice, set
$N_\Q:=N\otimes_Z \Q$ and $N_\R:= N\otimes_\Z\R$, and analogously define
$M_\Q$ and $M_\R$.

A \emph{rational polyhedral strictly convex cone} $\sigma \subset N_\R$ is a closed convex cone generated by finitely many vectors lying in $N$, and such that $\sigma \cap - \sigma = \{0\}$. Its dual cone $\check{\sigma} := \{ m \in M_\R, \, u(m) \ge 0  \text{ for all } u \in \sigma \}$ is a finitely generated semi-group. Thus $\sigma$ defines an affine variety $U_\s := \Spec \C[\check{\sigma}\cap M]$. 
The torus  $\bG_m^d = \Spec\C[M]$ is contained as a dense orbit in $U_\s$  and 
the action by $\bG_m^d$ on itself extends to $U_\s$, which makes $U_\s$ a toric variety. 
Conversely, any affine toric variety can be obtained in this way. 

If $\s$ is \emph{simplicial}, i.e.,  it is generated by vectors linearly independent over $\R$, then $U_\s$
 has at worst quotient singularities. The toric variety $U_\s$ is smooth if and only if $\s$ is simplicial and generated by $d$ vectors $e_1, \ldots , e_d$ forming a basis of  $N$ as an abelian group; such a $\sigma$ is said to be \emph{regular}.

A \emph{fan} $\Delta$ is a finite collection of rational polyhedral strictly convex cones in $N_\R$ such that each face of a cone in $\Delta$ belongs to $\D$ and the intersection of two cones in $\Delta$ is a face of both of them. A fan $\Delta$ determines a toric variety $X(\Delta)$, obtained by patching together the affine toric varieties $\{U_\sigma\}_{\s \in \Delta}$ along their intersections in a natural way.  
If all cones in $\Delta$ are simplicial then $\Delta$ is said to be
\emph{simplicial} and if all cones are regular, then $\Delta$ is said
to be \emph{regular}; $X(\D)$ is smooth if and only if $\D$ is regular. If $\bigcup_{\s\in\D}\sigma=N_\R$, then $\Delta$
is said to be \emph{complete}. The toric variety $X(\D)$ is compact if
and only if $\D$ is complete. Unless otherwise stated, we will assume that all fans in this paper are complete.

There is an one-to-one correspondence between cones of $\Delta$ of
dimension $k$ and orbits of the action of $\bG^d_m$ on $X(\Delta)$ of
codimension $k$.  We denote the closure of
the orbit associated with a cone $\s\in\Delta$
in $X(\D)$ by $X(\s)$. In particular, $1$-dimensional cones correspond to (irreducible) $\bG^d_m$-invariant divisors. 

A fan $\D'$ \emph{refines} another fan $\D$
if each cone in $\D'$ is included in a 
cone in $\D$.

\subsection{Equivariant (holomorphic) maps}
Given a group morphism $A: M\to M$, we will write $A$ also for the induced
linear maps $M_\Q\to M_\Q$ and $M_\R\to M_\R$. Morever, we let $\Atrans$
denote the 
dual map $N\to N$, as well as the dual linear maps $N_\Q\to N_\Q$ and
$N_\R\to N_\R$. 
It turns out to be convenient to use this notation rather than 
writing $\A$ for the map on $N$ and $\Atrans$ for the map on $M$. 

A \emph{map of fans} $\Atrans : (N,\D_2)\to  (N,\D_1)$ is a linear map
$\Atrans: N_\R \to N_\R$ that preserves $N$ and satisfies that 
the fan $\Atrans(\D_2):=\{\Atrans(\sigma): \sigma\in\D_2\}$ refines $\D_1$. 
If $\sigma_1\in\D_1$ and $\s_2\in\D_2$ satisfy that
$\Atrans(\s_2)\subseteq \s_1$, then the dual map $\A:M\to M$ maps $\check \s_1$ to $\check \s_2$ and induces a map $\A:\C[\check{\s}_1\cap M]\to \C[\check{\s}_2 \cap M]$, which, in turn, induces a map 
$\f_\A : U_{\s_2} \to U_{\s_1}$. These maps can be patched together to a holomorphic map  
$\f_\A : X(\D_2) \to X(\D_1)$ which  is \emph{equivariant} in the following sense: Denote by $\rho_\A : \bG^d_m \to \bG^d_m$ the natural group morphism induced by the ring morphism $\A: \C[M]\to \C[M]$.
Then for any $x\in X(\D_2)$, and any $g\in \bG_m^d$, one has
$ \f_\A ( g \cdot x) = \rho_\A(g) \cdot \f_\A(x)$. Conversely any equivariant holomorphic map $X(\D_2)\to X(\D_1)$ is determined by a map of fans $\Atrans : (N,\D_2)\to  (N,\D_1)$.

The map $\f_\A$ is dominant if and only $\det (\A)\neq 0$  and the topological degree of $\f_\A$ equals
$|\det (\A) |$.

%%%%%%%%%%%%%%%
\subsection{Universal cohomology of toric varieties}

Let $\Delta$ be a complete simplicial fan. Then $X(\Delta)$ has at worst quotient singularities, 
its cohomology groups $H^{j}(X(\Delta)) := H^{j}(X(\Delta),\R)$ with values in $\R$ vanish whenever $j$ is odd, and $H^*(X(\Delta))$ is generated 
as an algebra by the $\bG^d_m$-invariant divisors $[X(\s)]$, where
$\s$ runs over the $1$-dimensional cones in $\D$.

We let $\frD$  denote the set of all complete simplicial fans in $N$
and 
endow it with a partial ordering by imposing $\D \prec \D'$
if (and only if) $\D'$ refines $\D$. 
For any two fans $\D, \D' \in \frD$, one can find a third fan $\D''$ refining both; hence
$\frD$ is a directed set. 
Assume $\D \prec \D'$. Then the identity map on $N$ induces a map of fans
$\check \id_{\D',\D}:(N,\D')\to (N,\D)$, and thus yields 
a natural birational morphism $\pi:=\f_{\id_{\D',\D}}: X(\D')\to X(\D)$. This
map 
induces linear actions on cohomology, 
$\pi^* : H^*(X(\D)) \to H^*(X(\D'))$ and $\pi_* : H^*(X(\D')) \to
H^*(X(\D))$ that satisfy $\pi_*\pi^* = \id$; in particular, the map $\pi_*$ is surjective
and $\pi^*$ is injective.

The pushforward $\pi_*$ and pullback $\pi^*$ arrows make $\frD$ into
an inverse and directed set, respectively, 
and so the limits
\begin{equation*}
\Hproj:
=
\varprojlim_\frD H^*(X(\D))
\text{ and }
\Hinj:
=
\varinjlim_\frD H^*(X(\D))~
\end{equation*}
are well-defined infinite dimensional graded real vector spaces.  
We will refer to $\Hproj$ and $\Hinj$ as the \emph{universal  (inverse}
respectively, \emph{direct) cohomology} of toric varieties. 

In concrete terms, an element $\omega \in \Hproj$ is a collection 
of \emph{incarnations}
$\omega_\D \in H^*(X(\D))$ for each $\D\in \frD$, such that $\pi_*
(\omega_{\D'}) = \omega_{\D}$ if 
$\D \prec \D'$ and $\pi=\f_{\id_{\D',\D}}$.
An element $\omega \in \Hinj$ is determined by some class $\omega_\D
\in H^*(X(\D))$, and two classes  $\omega_{\D_i} \in H^*(X(\D_i))$, $i=1,2$ determine the same class in $\Hinj$
if and only if there exists a common refinement 
$\D' \succ
\D_i$ 
such that
$\pi_1^*(\omega_{\D_1}) =\pi_2^*(\omega_{\D_2})$ if $\pi_i=\phi_{\id_{\D',\Delta_i}}$.
Note that the map 
that sends $\omega_\D\in
H^*(X(\D))$ to the class it determines in $\Hinj$ is injective. 

We endow $\Hproj$ with its projective limit topology so that $\omega_j \to \omega$
if and only if $\omega_{j,\D} \to \omega_\D$ for each fan $\D\in\frD$. Then $\Hinj$ is dense in $\Hproj$.

Each cohomology space $H^*(X(\D))$ has a ring structure coming from
the intersection product  
 which respects the grading so that $\omega_\D \cdot\eta_\D\in H^{2(i+j)}(X(\D))$ if  $\omega_\D \in H^{2i}(X(\D))$ and
  $\eta_\D\in H^{2j}(X(\D))$. 
Given classes $\omega$ and  $\eta$
in $\Hinj$, 
pick $\Delta\in\frD$ so that they are determined 
by $\omega_\D$ and $\eta_\D$, respectively, and let $\omega \cdot \eta$ 
be the class in $\Hinj$ determined by $\omega_\D \cdot \eta_\D$. 
It is not difficult to check that this definition of $\omega\cdot\eta$ is independent of the choice of $\D$.
Hence, in this way, $\Hinj$ is endowed with
a natural structure of a graded $\R$-algebra.

More generally, given $\omega \in \Hproj$ and  $\eta\in\Hinj$, pick $\Delta$
such that $\eta$ is determined by $\eta_\D$ and let $\omega\cdot\eta$ be the
class in $\Hproj$ determined by $\omega_\D\cdot\eta_\D$. 
Again, this product is well-defined and independent of the choice of
$\D$ and so  
$\Hproj$ is a $\Hinj$-module. Note, however, that it is not possible to define 
a ring structure on $\Hproj$ that continuously extends the one on
$\Hinj$.

Since the intersection product $H^{2j}(X(\D))\times H^{2(d-j)}(X(\D))\to\R$ is a
perfect pairing for each $\D\in\frD$ by Poincar\'e duality,  the pairing
$\hproj^{2j} \times \hinj^{2(d-j)} \to \R$ is also perfect and thus
$\Hproj$ and $\Hinj$ are naturally dual one to the other.

%%%%%%%%%%%%%%%
\subsection{Toric line bundles}\label{sec:line}
The Picard group $\Pic X(\Delta)$ of a toric variety $X(\Delta)$ is generated by classes of $\bG^d_m$-invariant Cartier divisors.
These divisors can in turn be described in terms of functions on $N_\R$ as follows.
Let $\PL(\D)$ be the set of all continuous real-valued functions $h$ on $|\Delta|\subset N_\R$ 
that are \emph{piecewise linear with respect to $\Delta$}, i.e.,  such that  for any cone $\s\in \D$ there exists $m = m(\s) \in M$ with $h|_\s = m$. Given any $1$-dimensional cone $\sigma$ of $\D$, the associated \emph{primitive vector} is the first lattice point $u_i$ met along $\sigma=\R_{\geq 0} u_i$. Let $\D(1) = \{u_i\}\subset N$ be the set of primitive vectors of $1$-dimensional cones in $\Delta$. With $h\in\PL(\D)$, we associate the Cartier divisor $ D(h) := \sum h(u_i) \, X(\R_{\geq 0}  u_i)$. The map sending $h\in \PL(\D)$ to $\O(- D(h)) \in \Pic (X(\D))$ is surjective and the kernel is the space of linear
functions  $M\subset \PL(\D)$.  By taking the first Chern class, we get a linear map: 
\begin{equation}\label{chernmap}
\Theta_1: 
\PL(\D) \to H^2\left (X(\D)\right) , \,  h \mapsto \Theta_1(h) =  c_1 (\O(- D(h))~.
\end{equation}
When $X(\D)$ is smooth, the kernel of $\Theta_1$ is $M$ and the image is precisely $H^2(X(\D), \Z)$.  
Note that $\Theta_1$ extends by linearity to $\PL_\Q(\D):=\PL(\D)\otimes_\Z\Q$, corresponding to $\Q$-line bundles, with image $H^2(X(\D), \Q)$.

Let $\Atrans: (N,\D_2) \to (N, \D_1)$ be a map of fans, inducing a holomorphic map $\f_\A: X(\D_2) \to X(\D_1)$, and pick $h\in \PL(\D_1)$. Then the pullback $\f_\A^* D(h)$ is a well-defined Cartier divisor on $X(\D_2)$, equal to  
$D(h\circ \Atrans)$. It follows that 
\begin{equation}\label{teta}
\f_\A^* \Theta_1(h) = \Theta_1 (h \circ \Atrans). 
\end{equation}

There is a link between positivity properties of the classes in $H^2(X(\D))$ and convex geometry. A function $h\in\PL(\D)$ is said to be \emph{strictly convex (with
  respect to $\D$)} if it is convex and defined by different elements $h|_\sigma\in M$ for different $d$-dimensional cones $\sigma\in\D$. Recall that on a complete algebraic variety, a Cartier divisor $D$  is nef  if $D\cdot C \ge 0$ for any curve $C$.
The line bundle $\O(-D(h))$ over $X(\D)$ is nef if (and only if)
$h\in\PL(\D)$ is convex and it is ample if (and only if) $h$ is
strictly convex.

A function $h$ in $\PL(\D)$ determines a (non-empty) polyhedron 
$$
P(h) := \{ m \in M_\R , \,  m \le h  \} \subset M_\R~. 
$$
If $h$ is strictly convex with respect to some fan, then $P(h)$ is a compact \emph{lattice polytope} in $M_\R$, i.e.,  it is the convex hull of finitely many points
in the lattice $M$, and it has non-empty interior. 
 Conversely, if $P\subset M_\R$ is a lattice polytope, then 
the function 
\begin{equation}\label{hp}
h_P (u) := \sup \{ m(u), \, m \in P \} 
\end{equation}
is a piecewise linear convex function on $N_\R$. 
If $\D_P$ denotes the \emph{normal fan} of $P$ (see~\cite[Section 1.5]{F} for a definition)
then $h_P\in\PL(\D)$ precisely if $\D$ refines $\D_P$ and it is strictly convex
with respect to $\D_P$. If $h_P\in\PL(\D)$ then $\Delta$ is said to be \emph{compatible} with $P$.

If $\Atrans: (N,\D_2)\to (N,\D_1)$ is a map of fans, note that 
$(h\circ\Atrans)=A P(h)$.
%\begin{equation}\label{segt}
%(h\circ\Atrans)=A P(h). 
%\end{equation}
If $\D$ is compatible with $P$ and $Q$, then the intersection product 
\begin{equation}\label{snitt}
\O(-D(h_P))^k\cdot\O(-D(h_Q))^{d-k}= d!\Vol (P[k], Q[d-k]), 
\end{equation}
where $\Vol$ is the mixed volume as defined in Section \ref{mixedvol}, see \cite[p.~79]{Od}. 
\begin{eg}\label{pnex}
Take a basis $e_1, \ldots , e_d$ of $N$ with dual basis $e_1^*,\ldots, e_d^*$,  and set $e_0 := - \sum_1^d e_i$.
Let $\D$ be the unique fan whose $d$-dimensional cones are the $d+1$ cones 
$\s_i = \sum_{j \neq i} \R_{\geq 0} e_j$. Then $X(\D)$ is isomorphic to the projective space $\P^d$. 
Let $h$ be the unique function in $\PL(\D)$ that satisfies $h(e_0)=1$
and $h(e_i) =0$ if $i\ge1$; note that $h$ is strictly convex with
respect to $\D$. Then
$\O(- D(h)) = \O_{\P^d}(1)$ and moreover the polytope $P_{D(h)}$ is 
the standard simplex
$$
\Sigma_d := \{  u = \sum_1^d s_i e^*_i, \, s_i \le 0, \, \sum_1^d s_i \ge -1
\} \subset M_\R~.
$$
\end{eg}

%%%%%%%%%%%%%%%
\subsection{Piecewise polynomial functions}\label{sec:piece}

Higher cohomology classes of toric varieties in $\frD$ can be encoded in 
terms of (piecewise) polynomial functions on $N_\R$.

Given a fan $\D$, let $\PP(k,\D)$ be the set of 
\emph{piecewise polynomial functions (with respect to $\Delta$) of degree $k$}, i.e, %such 
continuous functions $h : N_\R \to \R$ 
such that for each cone $\s \in \D$, the restriction $h|_\s = \sum m_{i_1}\otimes \cdots  \otimes m_{i_k}$, $m_{i_j}\in M$, is a homogeneous polynomial of degree $k$.
Note that $\PP(1,\D)=\PL(\D)$. 
Moreover, note that $h\otimes h'\in\PP(k+k',\D)$ if $h\in\PP(k,\D)$ and
$h'\in\PP(k',\D)$, so that $\PP(\D) := \oplus_k \PP(k,\D)$ is
a graded ring.

From now on, assume $\D\in\frD$, and let $\sigma$ be a $j$-dimensional
cone in $\D$. Since $\s$ is
simplicial it is generated by exactly $j$ primitive vectors in
$\D(1)=\{u_1,\ldots, u_N\}$, say $u_1,\ldots, u_j$, and so $x\in\s=\sum_{i=1}^j\R_{\geq 0}u_i$ admits a unique
representation $x=\sum_{i=1}^j x_iu_i$ with $x_i\geq 0$. It follows that
$h|_\s$ has a unique expansion $h|_\s(x)= \sum a_I x_{i_1}\cdots x_{i_k}$, where the 
sum ranges over all $I=\{i_1,\ldots,
  i_k\}\subset\{1,\ldots, j\}$.

We can now define a linear map $\Theta_k: \PP(k,\D) \to
H^{2k}(X(\D))$ by 
$\Theta_k(h):=\sum a_IX_I$, 
where the sum ranges over all $I=\{i_1,\ldots, i_k\}\subset\{1,\ldots,
N\}$, $X_I$ is (the class of) the intersection product $X(\R_{\geq 0} u_{i_1})\cdot
\ldots \cdot X(\R_{\geq 0} u_{i_k})$, and $a_I$ is the coefficent of
$x_{i_1}\cdots x_{i_k}$ in $h|_\s$
defined above if $\s=\sum \R_{\geq 0}u_{i_\ell}$ is a cone in $\D$ and
$a_I=0$ otherwise. 

By patching the maps $\Theta_k$ together, we obtain a graded 
map
$\Theta: \PP(\D) \to H^*(X(\D)),  h\mapsto \sum_k \Theta_k(h_k)$, where 
$h_k$ is the $k$-th graded piece of $h$. 
Note that $\Theta$ is also a ring morphism since $\Theta (h\, h') =  \Theta (h) \, \Theta(h')$ for any $h,h'\in \PP(\D)$. If $X(\D)$ is smooth the image of $\Theta$ is $H^*(X(\D), \Z)$. As for $\Theta_1$ we can extend $\Theta$ to $\PP_\Q(\D):=\PP(\D)\otimes_\Z \Q$, with image $H^*(X(\D), \Q)$.

%%%%%%%%%%%%%%%
\subsection{Equivariant rational morphisms}\label{sec:morphisms}
Let $\A: M_\R \to M_\R$ be a linear map that preserves $M$. Take $\D,\D'\in\frD$ such that $\D'$ refines $\D$ and $\Atrans^{-1}(\sigma)$ is a union of cones in $\D'$ for each $\s\in\D$. 
Then $\Atrans:(N,\D')\to(N,\D)$ is a map of fans, inducing a holomorphic equivariant map $f_\A:X(\D')\to X(\D)$. 
Let $\pi:X(\D')\to X(\D)$ be the equivariant birational map induced by $\check \id_{\D',\D}:(N,\D')\to(N,\D)$, and let $\f_\A:= f_\A\circ \pi^{-1}$. Then $\f_\A : X(\D) \dto X(\D)$ is a rational map that is equivariant under the action of $\bG_m^d$. Conversely, any equivariant rational self-map on $X(\D)$ arises in this way. 
The map $\f_\A$ is holomorphic precisely if $\D\prec \Atrans(\D)$ and
it is dominant precisely if $\det (\A) \neq 0$.

Let $e_1,\ldots, e_d$ be a basis of $M$ and let $e_1^*,\ldots, e_d^*$ be the corresponding basis of $N$. Then $\A = \sum a_{ij}e_i\otimes e_j^*$, for some $a_{ij}\in \Z$. If  $x_1, \ldots , x_d$ are the induced coordinates on $\bG^d_m$, then $\f_\A$ restricted to  $\bG^d_m$ is the (holomorphic) monomial map $\f_\A(x_1, \ldots , x_d) =  (\prod x_j^{a_{j1}}, \ldots , \prod x_j^{a_{jd}})$.

Recall that a dominant holomorphic map $\f:X'\to X$ induces linear
actions on cohomology $\f^*:H^*(X)\to H^*(X')$ and $\f_*:H^*(X')\to
H^*(X)$. 
Assume that $\f_\A$ is dominant. Then we define the \emph{pushforward} $(\f_\A)_\bullet : H^*(X(\D))\to H^*(X(\D))$
as the composition $(\f_\A)_\bullet := (f_{\A})_* \circ \pi^*$, and the
\emph{pullback} $\f_{\A}^\bullet : H^*(X(\D))\to H^*(X(\D)$ as
$\f_\A^\bullet := \pi_* \circ f_\A^* $. It is readily verified that
$(\f_\A)_\bullet$ and $\f_\A^\bullet$ do not depend on the choice of
$\D'$. We insist on writing $(\f_{\A})_\bullet,\ \f_\A^\bullet $
instead of $(\f_{\A})_*,\ \f_\A^*$ since one does not have good
functorial properties, e.g. $(\f_B\circ\f_\A)_\bullet\neq
(\f_B)_\bullet \circ (\f_\A)_\bullet$ in general.

The linear map $\A$ also induces natural linear actions $\f_\A^*: \Hinj \to \Hinj$ and $(\f_\A)_* : \Hproj \to \Hproj$, defined as follows. 
Suppose $\eta$ is a class in $\Hinj$ determined by 
$\eta_{\D} \in H^*(X(\D))$. Pick $\frD\ni\D' \succ \D$ such that the map $f_\A : X(\D') \to X(\D)$ induced by $\Atrans$ is holomorphic, and define $\f_\A^*\eta$ to be the class in $\Hinj$ determined by $f_\A^*\eta_{\D} \in H^* (X(\D'))$. 
Next, suppose $\omega\in\Hproj$. The incarnation of $(\f_{\A})_* \omega$ in $H^*(X(\D))$ for a given $\D\in\frD$ is defined as $(\f_{\A})_*\omega_{\D} := (f_{\A})_*\omega_{\D'}$, where $\D'$ is choosen as above. It is not hard to check that $\f_\A^*$ and $(\f_\A)_*$ are independent of the choice of refinement $\D'$, and moreover that  
$(\f_{\A})_*$ is continuous on $\Hproj$,
$\f_{B\circ\A}^*=\f_\A^*\circ\f_B^*$, $(\f_{B\circ A})_*=(\f_B)_*\circ
(\f_\A)_*$, 
$(\f_{\A})_* \circ (\f_\A)^* = |\det (\A)|$, and 
$(\f_{\A})_* \omega \cdot \eta = \omega \cdot (\f_{\A})^* \eta$ 
for any two classes $\omega \in \hproj^{2j}$, $\eta\in \hinj^{2(d-j)}$.

Given $h\in\PL(\D)$, the first Chern class $\Theta_1(h)\in H^2(X(\D))$ determines a class in $\Hinj$, also denoted by $\Theta_1(h)$, that satisfies 
\begin{equation}\label{eq:pull-line}
\f_\A^* (\Theta_1(h)) = \Theta_1(h \circ \Atrans) \text{ in } \Hinj~, 
\end{equation}
which follows in light of \eqref{teta}.

%%%%%%%%%%%%%%%%%%%%%%%%%%%%%%%%%%%%%%%%%%%%%%%%%%%

\section{The polytope algebra}\label{sec:polytope}

%%%%%%%%%%%%%%%%%%

\subsection{Definition}\label{palgebra}
Given any finite collection  of convex sets $K_1,\ldots, K_s\subset M_\R$,  we let
$K_1+\cdots + K_s$ denote the 
\emph{Minkowski sum} $K_1 + \cdots + K_s:=\{x_1+\cdots +x_s \mid x_j\in K_j\}$. For any $r\in\R_{\geq 0}$, we also write $rK_j:=\{rx\mid x\in K_j\}$. A \emph{polytope in $M_\Q$} is the convex hull of finitely many 
points in $M_\Q$.

We now introduce the \emph{polytope algebra} $\Pi= \Pi(M_\R)$ 
which is a variant of the original construction
of P.~McMullen
~\cite{M}. It is the $\R$-algebra
with a generator $[P]$ for each polytope $P$ in 
$M_\Q$, with $[\emptyset]=:0$. The generators satisfy the relations 
$[P\cup Q]+[P\cap Q]=[P]+[Q]$ whenever $P\cup Q$ is convex, and $[P+t]=[P]$
for any $t\in M_\Q$. The multiplication in $\Pi$ is given by 
$[P]\cdot [Q]:=[P+Q]$, with 
multiplicative unit $1:=[\{0\}]$. 
The polytope algebra admits a grading 
$\Pi=\bigoplus_{k=0}^d\Pi_k$
such that $\Pi_k \cdot \Pi_l \subset \Pi_{k+l}$.  
The $k$-th graded piece $\Pi_k$ is the $\R$-vector space spanned by
all elements of the form $(\log [P])^k$, where 
$\log[P]:=\sum_{r=1}^d\frac{(-1)^{r+1}}{r}([P]-1)^r$ and $P$ runs over all
polytopes in $M_\Q$. The top-degree part $\Pi_d$ is one-dimensional,
and multiplication gives non-degenerate pairings
$\Pi_j\times\Pi_{d-j}\to\Pi_d$. Given $\a\in\Pi$, we will write
$\alpha_k$ for its homogeneous part of degree $k$.

The lattice $M$ determines a (canonical) volume element on $M_\R$, 
which we denote by $\Vol$. 
It is normalized by the convention $\Vol (P)  = 1$ for any
parallelogram $P = \{ \sum s_i e_i^*, \, 0\le s_i \le 1\}$ such that
$e_1^*,\ldots, e_n^*$ is a basis of the lattice $M$. In particular,
the volume of the standard simplex $\Vol (\Sigma_d )$ is $1/d!$. 
There is a canonical linear map $\Vol : \Pi \to \R$ defined by 
$\Vol ([P]) = \Vol (P)$. 
This map is zero on all pieces $\Pi_k$ for $k\le d-1$,
and it induces an isomorphism $\Vol : \Pi_d \mathop{\rightarrow}\limits^{\simeq} \R$.

\smallskip

Let $\A:M_\Q\to M_\Q$ be a linear map. 
Then $\A$ induces a linear map $\Pi\to\Pi$, defined by 
$[P]\mapsto [\A(P)]$; we shall denote it by
$\A_* : \Pi \to \Pi$. Note that $\A_*$ is actually a ring homomorphism since 
$\A_*([P]\cdot [Q])=[\A(P+Q)]=\A_*[P]\cdot \A_*[Q]$ for polytopes $P$ and
$Q$ in $M_\Q$, and $\A_*$ preserves the grading on $\Pi$ 
since $\A_*(\log[P])=\log[\A(P)]$. Also, it is clear that $(\A \circ \B)_*=\A_* \circ \B_*$ for any two linear maps  $\A$ and $\B$.

An important example is given by the homothety $\A=r\times\id$, $r\in \Q_{\geq 0}$; we denote the corresponding map on $\Pi$ by $\cD(r)$. 
Note that if $P$ is a polytope and $r\in\Z$, then 
$\cD(r)[P]=[P]^r$. 
It is proved in ~\cite[Lemma 20]{M} that 
\begin{equation}\label{homogeneity} 
\alpha\in\Pi_k  \text{ if and only if }
\cD(r)\alpha =r^k\alpha~,
\text{  for any } r \in \Q_{\geq 0} ~.
\end{equation}

If $\det (\A)\neq 0$, there is a well-defined 
\emph{pullback} map $\A^*:\Pi\to\Pi$ by $\A^*:=\cD\left(|\det (\A)|^{1/k}\right)\circ (\A^{-1})_*$ on
$\Pi_k$; in particular, 
\begin{equation*}
\A^*[P]=|\det (\A)|[\A^{-1}(P)]
\end{equation*} for any
polytope $P$. Moreover $(\A \circ \B)^*=\B^*\circ \A^*$ for any two
linear maps $A$ and $B$ with non-zero determinant. 
Beware that $\A^*$ is not a ring homomorphism on $\Pi$.
On the other hand, $\A^* ( \A_* (\alpha)) = |\det (\A)|\, \alpha$ for any $\a \in \Pi$.

%%%%%%%%%%%%%%%%%%%%%%%%

\subsection{Mixed volumes}\label{mixedvol}
Let $K_1,\ldots, K_s\subset M_\R$ be convex compact sets and pick $r_1,\ldots,
r_s\in\R_{\geq 0}$.  
A theorem by Minkowski and Steiner asserts that 
$\Vol(r_1K_1+\cdots +r_sK_s)$ is a homogeneous polynomial  of degree $d$ in the variables 
$r_1,\ldots, r_s$. In particular, there is a unique expansion: 
\begin{equation}\label{minkowski}
\Vol\left (r_1K_1+\cdots +r_sK_s\right )=
\sum_{k_1+\cdots +k_s=d} \binom{d}{k_1, \ldots, k_s}
\Vol\left (K_1[k_1],\ldots, K_s[k_s]\right )\,  r_1^{k_1}\cdots r_s^{k_s},
\end{equation}
the coefficients $\Vol(K_1[k_1],\ldots, K_s[k_s])\in \R$ are called
\emph{mixed volumes}. Here the notation $K_j[k_j]$ refers to the repetition of $K_j$ $k_j$
times. It is a fact that 
$\Vol(K_1[k_1],\ldots, K_s[k_s])$ is non-negative, multilinear
symmetric in the variables $K_j$, and increasing in each variable,
meaning that  
\begin{equation}\label{monotone}
\Vol\left (K_1[k_1], K_2[k_2], \ldots, K_s[k_s]\right ) \leq \Vol\left (K'_1[k_1],
K_2[k_2], \ldots, K_s[k_s]\right ) \text{ whenever } K_1\subseteq K_1'.
\end{equation}
Note that $\Vol(K_1[d], K_2[0], \ldots, K_s[0])=\Vol(K_1)$. 
There is in general no simple geometric description of mixed volumes,
unless the $K_j$ has some symmetries, cf. Section~\ref{sec:diag} and \eqref{eq:cauchy-crofton} below.

Since $\Vol(K)$ is invariant under translation of $K$, 
\begin{equation}\label{trans}
\Vol \left ((K_1+t)[k_1], K_2[k_2],\ldots, K_s[k_s]\right ) = \Vol \left (K_1[k_1], K_2[k_2],\ldots, K_s[k_s]\right )~,
\end{equation}
for any $t\in M_\R$.
Moreover, $\Vol(K_1[k_1],\ldots, K_s[k_s])\in \R$ is additive in the sense that 
\begin{multline*}
\Vol \left ((K_1\cup K_1')[k_1], K_2[k_2],\ldots, K_s[k_s]\right ) + \Vol \left ((K_1\cap K_1') [k_1], K_2[k_2],\ldots, K_s[k_s]\right ) 
=\\ \Vol \left (K_1[k_1], K_2[k_2],\ldots, K_s[k_s]\right ) + \Vol \left (K_1'[k_1], K_2[k_2],\ldots, K_s[k_s]\right );
\end{multline*}
as soon as $K_1\cup K_1'$ is convex.
It follows that the mixed volumes extend to the polytope algebra $\Pi$ 
as multilinear functionals:
$$
\Pi^s\ni (\a_1, \ldots , \a_s) \mapsto \Vol \left (\a_1[k_1],\ldots, \a_s[k_s]\right ) \in \R~,
$$
so that, in particular, $\Vol\left ([P_1][k_1],\ldots,
[P_s][k_s]\right )=\Vol\left (P_1[k_1],\ldots, P_s[k_s]\right)$. 
Equation~\eqref{minkowski} translates into 
\begin{equation}\label{eq:mink}
\Vol\left (\cD(r_1)\alpha_1\cdot \ldots \cdot \cD(r_s)\alpha_s\right ) = 
\sum_{k_1+\ldots +k_s=d} {\binom{d}{k_1,\ldots, k_s}}\Vol\left (\alpha_1 [k_1], \ldots,
\alpha_s[k_s]\right )\,  r_1^{k_1}\cdots r_s^{k_s}~,
\end{equation}
which holds for $r_1,\ldots, r_s\in\Q_{\geq 0}$. 
Note that \eqref{eq:mink} implies the following homogeneity: 
\begin{equation}\label{eq:mink-hom}
\Vol\left (\cD(r_1)\alpha_1[k_1],  \ldots ,  \cD(r_s)\alpha_s[k_s]\right ) = 
\Vol\left (\alpha_1 [k_1], \ldots,
\alpha_s[k_s]\right )\,  r_1^{k_1}\cdots r_s^{k_s}~.
\end{equation} 

\begin{lem}\label{lem:zero}
Let $\alpha_1, \ldots, \a_s$ be homogeneous elements in the polytope
algebra of degrees $\ell_1,\ldots, \ell_s$, respectively. 
Then  $\Vol\left (\alpha_1 [k_1], \ldots,
\alpha_s[k_s]\right )=0$ unless $\ell_j = k_j$ for all $j$, in which case it is equal to 
${\binom{d}{\ell_1,\ldots, \ell_s}}^{-1}\, \Vol \left (\a_1\cdot \ldots \cdot \a_s\right )$.
\end{lem}
\begin{proof}
By~ \eqref{homogeneity}, and the linearity of $\Vol:\Pi\to\R$,
\begin{equation*}
\Vol\left (\cD(r_1)\alpha_1\cdot \ldots \cdot \cD(r_s)\alpha_s\right )
= \Vol \left (r_1^{\ell_1} \a_1\cdot \ldots \cdot r_s^{\ell_s} \a_s\right ) =
 \Vol \left (\a_1\cdot \ldots \cdot \a_s\right )\, r_1^{\ell_1}\cdots r_s^{\ell_s};
\end{equation*}
in particular, the only non-vanishing mixed volume is $\Vol\left (\a_1[\ell_1],\ldots, \a_s[\ell_s]\right )={\binom{d}{\ell_1,\ldots, \ell_s}}^{-1}\, \Vol \left (\a_1\cdot \ldots \cdot  \a_s\right )$.
\end{proof}
Lemma ~\ref{lem:zero} implies that if $\a_1,\ldots,
\a_s\in\Pi$, and $\a_{j,\ell}$ denotes the $\ell$-th graded part of
$\alpha_j$, then 
\begin{equation}\label{madison}
\Vol\left (\a_1[k_1],\ldots, \a_s[k_s]\right )=
{\binom{d}{k_1,\ldots, k_s}}^{-1}\, \Vol \left (\a_{1,k_1}\cdot \ldots
\cdot \a_{s,k_s}\right )
\end{equation}

\begin{lem}\label{lem:dual}
Let $\A: M_\Q\to M_\Q$ be a linear map such that $\det(\A)\neq0$. Then 
$$
\Vol \left (\A^*\alpha_1 [k], \alpha_2[d-k]\right ) = \Vol \left (\alpha_1 [k], \A_*
\alpha_2[d-k]\right ), $$
for any two elements $\a_1, \a_2\in \Pi$.
\end{lem}

\begin{proof}
By multlinearity, we may assume that $\a_i =[P_i]$ for some polytopes
$P_i$.
Note that for $r_i\in\Q_{\geq 0}$,
\begin{equation*}
\Vol \left ( r_1 \A^{-1}(P_1) + r_2 P_2\right )=
\sum_{k} {\binom{d}{k}}
\Vol\left (\A^{-1} (P_1)[k], P_2[d-k]\right )\,  r_1^{k}r_2^{d-k},
\end{equation*}
and 
\begin{multline*}
|\det(\A)|\, 
\Vol \left ( r_1 \A^{-1}(P_1) + r_2 P_2\right ) = 
\Vol \left ( \A (r_1 \A^{-1}(P_1) + r_2 P_2)\right ) = \\
\Vol \left ( r_1P_1 + r_2 \A(P_2)\right )  =
\sum_{k} {\binom{d}{k}}
\Vol\left (P_1[k], \A(P_2) [d-k]\right )\,  r_1^{k}r_2^{d-k}.
\end{multline*}
By identification of the coefficients of $r_1^kr_2^{d-k}$ 
we get
\begin{equation}\label{ident}
|\det(\A)|\, \Vol\left (\A^{-1} P_1[k], P_2[d-k]\right ) = 
\Vol\left (P_1[k], \A(P_2) [d-k]\right ).
\end{equation}
The right hand side of \eqref{ident} is precisely
$\Vol \left ([P_1] [k], \A_* [P_2][d-k]\right )$, and in light of Lemma
~\ref{lem:zero}, the left hand side is
equal to 
\begin{multline*}
|\det(\A)|\, \Vol\left ([\A^{-1} (P_1)]_k[k], P_2[d-k]\right ) = 
\Vol\left[\cD\left( |\det(\A)|^{1/k}\right)[\A^{-1} (P_1)]_k [k], P_2[d-k]\right]=\\
\Vol\left (\A^* [ P_1]_k [k], P_2[d-k]\right ) = \Vol\left ( \A^* [ P_1] [k],
P_2[d-k]\right ), 
\end{multline*} 
which concludes the proof. Here we have used \eqref{eq:mink-hom}, the
definition of $\A^*$, and Lemma ~\ref{lem:zero}  for the
first, second, and last equalities, respectively. 
\end{proof}

 %%%%%%%%%%%%%%%%%%
 
\subsection{Currents on the polytope algebra}\label{sec:current} 
In order to simplify computations and relate the polytope algebra to the 
universal cohomology of toric varieties it is 
convenient to introduce the following terminology. A \emph{current} is a linear form on the polytope algebra. 
We denote the space of currents on $\Pi$ by $\CC$ and endow it with
the topology of pointwise convergence. Moreover, we write 
$\langle T , \b \rangle\in \R$ for the action of $T\in\mathcal C$ on $\b\in \Pi$.

A current $T\in\mathcal C$ is said to be of \emph{degree $k$} if
$T|_{\Pi_j}=0$ for $j\neq d-k$. Let $\mathcal C_k$ denote the subspace
of $\mathcal C$ of currents of degree $k$.  Note that $T\in\mathcal C$ admits
a unique decomposition $T=\sum T_k$, where $T_k\in\mathcal C_k$. (In
fact, $T_k$ is the trivial extension to $\Pi$ of the restriction of the linear
form $T$ to $\Pi_k$.)

Any invertible linear map $\A:M_\Q\to M_\Q$ induces actions on
$\mathcal C$, dual to the pullback and pushforward on $\Pi$, 
defined by 
$\langle \A_* T , \b \rangle := \langle T ,\A^* \b \rangle$ and
$\langle \A^* T , \b \rangle := \langle T ,\A_* \b \rangle$ for
$T\in\mathcal C$ and $\b\in\Pi$. 
It is not difficult to see that $T\in\CC$ is homogeneous of degree
$k$ if and only if 
$\cD(r)^* T = r^k T$ for any $r\in\Q_{\geq 0}$.

\smallskip

Let us describe some important examples of currents. 
\begin{eg}\label{poly-embedding}
Pick  $\a\in \Pi$, and let $T_\a$ be the current defined by 
$T_\a(\b) := \Vol (\a\cdot \b)$ for $\b\in\Pi$. The map $\a \mapsto
T_\a$ gives a linear injective map $\Pi\to\CC$ that sends $\Pi_k$ to
currents of degree $k$. 

In general, for $\a = \sum \a_k\in \Pi$, with $\a_k\in\Pi_k$, 
$T_{\a_k}(\b) = {d \choose n}\Vol ( \a[k], \b[k-d])$, 
which follows immediately from Lemma~\ref{lem:zero}. Moreover, by 
Lemma~\ref{lem:dual}, the actions of an invertible linear map $\A: M_\Q\to M_\Q$ on $\Pi$ and $\CC$ are compatible so that
$T_{\A^*\a} = \A^*T_{\a}$ and $T_{\A_*\a} = \A_*T_{\a}$ for any class $\a\in \Pi$.
\end{eg}

\begin{eg}\label{valuationex}
Given a vector space $V$, we define the \emph{convex body algebra}  $\K(V)$ as the polytope algebra, but with a generator $[K]$ for each compact 
convex set $K\subset V$ and with the relations $[K+t] = [K]$ for any $t\in V$, and $[K\cup L]+[K\cap L]=[K]+[L]$ whenever $K\cup L$ is convex. 
A \emph{(continuous translation-invariant) valuation} is a linear map on
$\K(V)$ that is continuous for the Hausdorff metric on compact sets, see e.g. ~\cite[Sect.~3.4]{S}. 
Let $\Val(V)$ denote the space of valuations on the space of convex bodies in $V$. 
Restricting the action of valuations on $M_\R$ to $\Pi$ gives 
an injective morphism $0\to\mathrm{Val}(M_\R)\to\CC$.
The construction of the current $T_\alpha$ in Example
~\ref{poly-embedding} can be extended to $\alpha,\beta\in\K(M_\R)$, and 
the mapping $\alpha\mapsto
T_\alpha$ embeds $\K(M_\R)$ into $\Val (M_\R)$. Thus
$\Pi \subset \K(M_\R)  \subset \mathrm{Val}(M_\R) \subset \CC~.$
\end{eg}

\begin{eg}\label{proj-val}
Endow $M_\R$ with an arbitrary euclidean metric $g$.
Let $H$ be a linear subspace of $M_\R$ of codimension $k$, and
$\Vol_H$ be the volume element on $H$ induced by $g$. Consider any linear projection 
$p: M_\R \to H$ onto $H$. Since $p(P\cap Q)=p(P)\cap p(Q)$ whenever
$P\cup Q$ is convex, $p$ can be extended to a function $\Pi\to\Pi$,
defined by $p[P]:=[p(P)]$, and the  linear map $\a\mapsto \Vol_{H} (p(\a))$ is a valuation of degree $k$ that we shall denote by $[H,p]$.
\end{eg}

%%%%%%%%%%

\subsection{Relations to the inverse cohomology of toric varieties}\label{sec:morphism}

Each polytope $P$ in $M_\Q$ determines a class $\ch(P)$ in $\Hinj$,
defined as follows:  Let $h_P$ be defined as in \eqref{hp} and 
choose $\D\in\frD$ so that $h_P\in\PL_\Q(\D)$. 
Now $\ch(P)$ is determined by 
the Chern character of the 
associated $\Q$-line
bundle
\begin{equation}\label{eq:defon pol}
(\ch  (P))_\D := \sum_{k=0}^d \frac1{k!} \Theta_1(h_P)^k  = \Theta\left( \sum_{k=0}^d \frac1{k!} h_P^k \right) = \Theta (\exp (h_P)) \in H^*(X(\D))~,
\end{equation}
where $\Theta_1$ and $\Theta$ are as in Sections ~\ref{sec:line} and
~\ref{sec:piece}, respectively.

The Chern character induces a linear map from the vector space $\oplus_P \R [P]$ to
$\Hinj$, defined by $\ch(\sum t_j[P_j])=\sum t_j \ch(P_j)$. 
We claim, in fact, $\ch$ is a well-defined ring homomorphism from
$\Pi$ to $\Hinj$. 
To see this, first note that $h_{P+t}=h_P+t$ for $t\in M_\Q$. 
It follows that $D(h_{P+t})$ and
$D(h_P)$ are linearly equivalent, see Section ~\ref{sec:line}. In
particular, $\Theta_1(h_{P+t})=\Theta_1(h_P)$, which implies $\ch(P+t)=\ch(P)$. 
Next, if $P$ and $Q$ are polytopes in $M_\Q$ such that $P\cup Q$ is convex, then
$h_{P\cup Q} = \max \{ h_P, h_Q\}$ and $h_{P\cap Q} = \min \{ h_P,
h_Q\}$. Thus 
$h_{P\cup Q}^k + h_{P\cap Q}^k = h_P^k + h_Q^k$ 
for any $k\ge0$, and by linearity of $\Theta$, 
$\ch  ([P\cap Q]+ [P\cup Q]) = \ch ([P]+[Q])$. 
Thus $\ch:[P]\mapsto \ch(P)$ is well-defined. 
Next, note that if $P$ and $Q$ are polytopes in $M_\Q$, then $h_{P+Q}=h_P+h_Q$. Hence
 $\ch([P]\cdot [Q])=\ch([P+Q])=\ch(P + Q) = \ch (P) \ch (Q) =\ch([P])\ch([Q])$; and so the claim is proved.

Let $\deg: \Hinj \to \R$ be the linear \emph{degree} map 
that is $0$ on $\hinj^k$
for $k < d$ and sends the class determined by a point in $X(\D)$ to $1$. 
The following theorem is due to
Fulton-Sturmfels~\cite[Sect.~5]{FS} and Brion~\cite[Sect.~5]{B}. 
\begin{thm}\label{thm:equiv} 
The  Chern character map $\ch:[P]\mapsto \ch(P)$ is an isomorphism 
of graded algebras 
$\ch : \Pi \to \Hinj$. 
It holds that $\deg  (\ch (\a)) = \Vol (\a)$ for $\alpha\in\Pi$.

By duality, we get a continuous isomorphism  
$\tch :  \Hproj \to \CC$, 
defined by 
$\langle \tch  (\omega) , \b \rangle: = \omega \cdot \ch  (\b)$ for $\omega \in\Hproj$ and $\b \in \Pi$.  
\end{thm}

\begin{prop}\label{chernprop}
Let $\A: M_\Q \to M_\Q$ be a linear map with $\det(\A) \neq 0$. 
Then 
\begin{equation}\label{eq:com1}
\ch  (\A_* \a ) = \f_\A^* \ch (\a)    
\end{equation}
for any $\alpha\in\Pi$. 
Similarly, 
for any $\eta\in \Hproj$, 
\begin{equation}\label{eq:com2}
\tch  ((\f_{\A})_* \eta ) = \A^* (\tch (\eta)). 
\end{equation}
\end{prop}

\begin{proof}
By linearity we may assume that 
$\a = [P]$ for some polytope $P$ in $M_\Q$. By definition, $\A_*[P] = [\A(P)]$ and $\ch  ([\A(P)])$ is the class in $\Hinj$ determined by $\Theta(\exp (h_{\A(P)}))$. 
Now 
$$
h_{\A(P)} 
= \sup\{ m, \, m \in \A(P)\} 
= \sup\{ m\circ \Atrans, \, m \in P\} 
= h_P \circ \Atrans. 
$$

In light of ~\eqref{eq:pull-line} and ~\eqref{eq:defon pol}, it follows that $\ch  ([\A(P)])$ is the pullback under $\f_\A$ of the class determined by $\Theta(\exp(h_P))$, i.e., $\ch  ([\A(P)]) = 
\f_\A^* \ch  ([P])$.

Now \eqref{eq:com2} follows from \eqref{eq:com1} by duality. Indeed, for $\eta\in \Hproj$ and $\a\in \Pi$, we have 
\begin{equation*}
\langle 
\tch  \left ((\f_{\A})_*\eta\right ) , \a
\rangle
=
(\f_{\A})_*\eta \cdot  \ch  \a
=
\eta \cdot \f_{\A}^* \ch  \a
= 
\eta \cdot  \ch (\A_* \a)
=
\langle 
\tch  (\eta),  \A_* \a
\rangle
=
\langle 
\A^*\tch  (\eta),   \a
\rangle.
\end{equation*} 
Here we have used the definition of $\tch$ for the first and fourth equality and \eqref{eq:com1} for the third equality. 
Moreover, the second and the last equality follow by Sections ~\ref{sec:morphisms} and ~\ref{sec:current}, respectively. 
\end{proof}

 %%%

\section{Dynamical degrees on toric varieties}\label{sec:dynamical} 
Let $\D$ be a complete simplicial fan and let $h$ be a strictly convex
piecewise linear function with respect to $\D$. Furthermore, let $\A:M\to M$ be
a group morphism and let 
$\f:=\f_\A:X(\D)\dto X(\D)$ be the corresponding rational equivariant map. 
The \emph{$k$-th degree of $\f$ with respect to} the
ample divisor $D:=D(h)$ is defined as 
$$
\deg_{D,k} (\f) :=  \f^\bullet D^k \cdot D^{d-k} \in \R_{\geq 0}.
$$
If $X(\D)=\P^d$ and $\O(D)=\O_{\P^d}(1)$, then $\deg_{D,k}(\phi)$ coincides with the
\emph{$k$-th degree of $\f$} $\deg_k(\f)$ as defined in the
introduction (Section 1).

The following result is a key ingredient in the proofs of Theorems ~A and ~D. 
\begin{prop}\label{volume}
Let $\D$ be a complete simplicial fan and let $D$ be an ample $\mathbb
G_m^d$-invariant divisor
on $X(\D)$ with corresponding polytope $P_D$. Moreover, let $A:M\to M$ be a group morphism with $\det
(A)\neq 0$, and let $\f_A:X(\D)\dashrightarrow X(\D)$ be the
corresponding equivariant rational map. 
Then 
\begin{equation}\label{eq:interpret}
\deg_{D,k} (\f_\A)  = d! \Vol ( \A( P_D) [k], P_D[d-k])~.
\end{equation}
\end{prop}

Recall from Example ~\ref{pnex} that if $X(\D)=\P^d$ and $D=\O_{\P^d}(1)$, then
$P_D$ is the standard simplex $\Sigma_d$. 
In this case \eqref{eq:interpret} reads
\begin{equation*}
\deg_{k} (\f_\A)  = d! \Vol ( \A( \Sigma_d) [k], \Sigma_d[d-k])~.
\end{equation*}

\begin{proof}
Pick a complete simplicial fan $\D'\succ\D$ such that the dual $\Atrans:N\to N$ of $A$ is a map of fans $\Atrans:(N,\D')\to (N,\D)$, let $f_\A: X(\D')\to X(\D)$ be the corresponding
equivariant map, and let $\pi: X(\D')\to X(\D)$ be the map induced
by $\check \id_{\Delta',\Delta}$.

Recall from Section ~\ref{sec:morphisms}
that then $\f_\A^\bullet = \pi_* \circ f_\A^*$. Hence
\begin{equation}\label{hoho}
\deg_{D,k} (\f_\A)= \pi_* \circ f_\A^* (D)^k \cdot D^{d-k} = f_\A^* (D)^k \cdot \pi^* (D^{d-k})
= (f_\A^* D)^k \cdot  (\pi^* D)^{d-k}.
\end{equation}
By Section \ref{sec:line}, and in particular \eqref{snitt}, the right hand side of \eqref{hoho} equals 
$$
 d! \Vol (\A(P_D)[k], P_D[d-k]),
$$
which concludes the proof. 
%Now $c_1(\O(D)) \in H^2(X(\D))$ and $c_1(\O(\pi^* D))\in H^2(X(\D'))$ determine the same class in $\Hinj$, which we denote by $[D]$, 
%and  $f_\A^* c_1(\O(D))\in H^2(X(\D'))$ determines the class
%$\f_\A^*[D]\in\Hinj$. Thus, in light of Sections
%~\ref{sec:morphisms} and ~\ref{sec:morphism}, 
%$\deg_{D,k} (\f_\A)$ is the degree of the 
%intersection product 
%of $(\f_\A^*[D])^k\in \hinj^k$, and  $[D]^{d-k}\in \hinj^{d-k}$.
%Note that $D=\Theta_1(h_{P_D})\in H^*(X(\D))$; thus $[D]^k =
%k!(\ch[P_D])_k\in\Hinj$. Moreover, by (the proof of) Proposition
%~\ref{chernprop}, $f_\A^*D=\Theta_1(h_{\A(P)})\in H^*(X(\D'))$, and 
%so $(\f_\A^*[D])^{d-k}=(d-k)! (\ch[\A(P)])_{d-k}\in\Hinj$. 
%Using this we get 
%\begin{multline*}
%\deg_{D,k}(\f_\A)=
%k!(d-k)!\deg((\ch  [\A(P_D)])_k \cdot (\ch  [P_D])_{d-k})=
%\\
%k!(d-k)!\Vol( [\A(P_D)]_k \cdot [P_D]_{k-d}) = d! \Vol (\A(P_D)[k], P_D[d-k])
%\end{multline*}
%Here we have used Theorem ~\ref{thm:equiv} for the second equality and
%\eqref{madison} for the third equality.
\end{proof}

Let us collect some  basic properties of $k$-th degrees.
These results are well-known and valid for arbitrary rational maps, see~\cite{DS}. However the case of toric maps is particularly simple. 
%Since it illustrates the power of the identification of  cohomology classes with elements of the polytope algebra, we shall provide full proofs of these statements.

\begin{prop}\label{prop:compare}
Let $\D_1$ and $\D_2$ be complete simplicial fans, and let $D_1$ and $D_2$ be
ample $\mathbb G_m^d$-invariant divisors on $X(\D_1)$ and $X(\D_2)$,  respectively. 
Then there exists a constant $C$ such that for any group morphism
$\A : M \to M$, one has 
\begin{equation}\label{chicago}
C^{-1} \, \deg_{D_2,k} (\f_{\A,2})
\le
\deg_{D_1,k} (\f_{\A,1})
\le 
C \, \deg_{D_2,k} (\f_{\A,2})
\end{equation}
where $\f_{\A,i}: X(\D_i) \dto X(\D_i)$, $i=1,2$ denote the respective induced maps.
\end{prop}

\begin{proof}
We claim that there is a constant $C$ such that
\begin{equation*}
C^{-1}\Vol(\A(P_{D_2})[k],P_{D_2}[d-k])\leq
\Vol(\A(P_{D_1})[k],P_{D_1}[d-k])\leq C\Vol(\A(P_{D_2})[k],P_{D_2}[d-k]).
\end{equation*}
Then \eqref{chicago} follows immediately from Proposition ~\ref{volume}. 

Since $D_1$ and $D_2$ are ample, and since $\Vol$ is translation
invariant, \eqref{trans}, in order to prove the claim we may assume that $P_{D_1}$ and $P_{D_2}$
contain the origin in $M_\R$ in their interior. Then for some $C_0$
large enough, 
$C_0^{-1}\, P_{D_2}\subset P_{D_1} \subset C_0\, P_{D_2}$. It
follows that the claim holds for $C=C_0^d$ since $\Vol$ is multilinear and
monotone, \eqref{monotone}.
\end{proof}

\begin{prop}\label{prop:submult}
Let $\D$ be a complete simplicial fan and let $D$ be an ample $\mathbb
G_m^d$-invariant divisor on $X(\D)$.
Then there exists a constant $C$ such that for any group morphisms
$\A_1, \A_2 : M \to M$, 
one has $$
\deg_{D,k} (\f_{\A_1}\circ \f_{\A_2})
\le 
C \, \deg_{D,k} (\f_{\A_1}) \, 
 \deg_{D,k} (\f_{\A_2}).$$
 \end{prop}

\begin{proof}
Given a rational map $f:\P^d\dto\P^d$, we
denote by $C(f)$ the set of  points  $p\in\P^d$ that  are either  indeterminate or critical for $f$, and by $PC(f) := \pi_2 \pi_1^{-1} (C(f))$ where $\pi_1, \pi_2$ denote the two projections of the graph of $f$ onto $\P^d$. This defines two proper algebraic subsets of $\P^d$.

\label{def-pullback}If $Z$ is  a variety of pure codimension $k$ in $\P^d$, we
denote by $f^{-1}(Z)$ the closure in $\P^d$ of  $f^{-1} (Z \cap PC(f))$.
Note that  by construction,  $f^{-1}(Z)$ is of codimension $k$ (or empty). We have the general inequality
$\deg (f^{-1} (Z)) \le \deg_k(f) \, \deg(Z)$, and for a generic
choice of $Z$, $\deg (f^{-1} (Z)) = \deg_k(f) \, \deg(Z)$.  
In particular, if $L$ is a generic linear subspace of $\P^d$ of codimension $k$,
then $\deg_k(f)=\deg(f^{-1}(L))$. 

We always have $\f_{\A_1}^{-1} ( \f_{\A_2}^{-1} (L))
=  (\f_{\A_1}\circ \f_{\A_2})^{-1} (L)$ outside $W: = C(\f_{\A_2}) \cup \f_{\A_1}^{-1} C (\f_{\A_1})$.
For a generic choice of $L$, the closure of  $(\f_{\A_1}\circ \f_{\A_2})^{-1} (L) \cap W$
is equal to $(\f_{\A_1}\circ \f_{\A_2})^{-1} (L)$. 
Whence
\begin{equation}\label{pellegrino}
\deg ((\f_{\A_1}\circ \f_{\A_2})^{-1}(L)) 
= \deg (\f_{\A_2}^{-1} ( \f_{\A_1}^{-1}(L))) \le \deg_k (\f_{\A_2}) \, \deg (\f_{\A_1}^{-1}(L))
\end{equation}
Since $L$ is generic the left hand side of \eqref{pellegrino} equal $\deg_k(\f_{\A_1}\circ \f_{\A_2})$ and the right hand side equals $\deg_k (\f_{\A_2}) \, \deg_k (\f_{\A_1})$. Thus $\deg_k(\f_{\A_1}\circ \f_{\A_2}) \le \deg_k(\f_{\A_1})\, \deg_k(\f_{\A_2})$, and applying Proposition ~\ref{prop:compare} to $D_1 = \O_{\P^d}(1)$ and $D_2 = D$, we get 
\begin{multline*} \deg_{D,k} (\f_{\A_1}\circ \f_{\A_2})
\le  C \, \deg_{k} (\f_{\A_1}\circ \f_{\A_2})
\le C\,  \deg_{k} (\f_{\A_1}) \,
 \deg_{k} (\f_{\A_2})\le \\
 C^3\,  \deg_{D,k} (\f_{\A_1}) \, 
 \deg_{D,k} (\f_{\A_2}), 
 \end{multline*}
 which concludes the proof.
\end{proof}

Pick a group morphism $\A: M\to M$, a fan $\D\in\frD$, and an ample
$\mathbb G_m^d$-invariant divisor $D$ on $X(\D)$. Then Proposition ~\ref{prop:submult} implies $$C\, \deg_{D,k} (\f_{\A}^{n+m}) \le (C\, \deg_{D,k} (\f_{\A}^{n})) \, (C\,  \deg_{D, k} (\f_{\A}^{m}))~.$$
Since the sequence $\{C\, \deg_{D,k}(\f_\A^n)\}_n$ is
sub-multiplicative, and $C^{1/n} \to 1$, we can define the
\emph{$k$-th dynamical degree of $\f_\A$ with respect to $D$}, 
$$
\l_{D,k}(\f_\A) := \lim_n \deg_{D,k} (\f_{\A}^{n})^{1/n}.  
$$

Assume $\D_1, \D_2 \in \frD$ and that $D_1, D_2$ are ample $\mathbb
G_m^d$-invariant divisors on  
$X(\D_1)$ and $X(\D_2)$, respectivly. Let $A:M\to M$ be a group
morphism, and let $\f_{\A,i}: X(\D_i)\dto X(\D_i)$, $i=1,2$ denote the induced equivariant
morphisms. Then Proposition ~\ref{prop:compare} implies that $\l_{D_1, k}(\f_{\A_1})=\l_{D_2, k}(\f_{\A_2})$. 
We shall write $\l_k(\f_\A)$ for the \emph{$k$-th dynamical degree of
  $\f_\A$} (computed in any toric model, and with respect to any ample
divisor).

For the record, we mention the following properties of the dynamical
degrees. 
Proposition ~\ref{volume} applied to $P_D$ yields that $\deg_{D,0}(\f_\A^n)  =  d! \Vol ( P_D)$ for all $n$ and $\deg_{D,d}(\f_\A^n)  = d! \Vol (\A^n(P_D))= d! |\det(\A)|^n \,
\Vol(P_D)$; therefore $\l_0(\f_A) =1$ and $\l_d (\f_A)= |\det(\A)|$. 
Moreover $\l_1(\f_\A) = \rho(\A)$, where $\rho(\A)$ is the spectral radius of $A$, i.e., the largest modulus of an eigenvalue of $A$; a proof is given in~\cite[Sect.~6]{HP}. 
\begin{prop}\label{specfall}
Let $\A: M\to M$ be a group morphism 
and $\D\in\frD$, and denote by $\f_\A :X(\D)\to X(\D)$ the induced equivariant morphism. Then, for any $0\le k,l \le d$,
\begin{equation}\label{cremona}
\l_{k+l}(\f_\A)\le \l_k(\f_\A)\, \l_l(\f_\A).
\end{equation} 
\end{prop}

\begin{proof}
By the Aleksandrov--Fenchel inequality, see (6.4.5) on p.~334 of~\cite{S},  
$$
\Vol ( \A (P_D)[k+l], P_D[d-k-l])\, \Vol(P_D)
\le
\Vol ( \A(P_D)[k], P_D[d-k])
\, 
\Vol ( \A(P_D)[l], P_D[d-l])
$$
which, in light of Proposition ~\ref{volume} implies \eqref{cremona}.
\end{proof}
Note that Proposition ~\ref{specfall} also immediately follows from Corollary B, taking it for granted.

%%%%%%%%%%

\section{Degree growth - Proof of Theorem A}\label{sec:pfA}

The proof of Theorem A can be reduced to controlling the growth of mixed volumes of convex bodies under the action of a linear map.  
Indeed, let $\A:M \to M$ be a group morphism and let $D$ be a divisor on a toric variety $X(\D)$, where $\D$ is a complete simplicial fan. Then,  
 by Proposition ~\ref{volume}, $\deg_{D,k}(\f_\A)=d!\Vol(\A(P_D)[k],
 P_D[d-k])$. In particular, 
$\deg_{k}(\f_\A)=d!\Vol(\A(\Sigma_d)[k], \Sigma_d[d-k])$, where $\Sigma_d$ is the standard simplex, see Example ~\ref{pnex}. Now  
Theorem A follows immediately from the following result. 
\begin{thm}\label{thm:cvx1}
Let $\A : M \to M$  be  a group morphism such that $\det(\A) \neq 0$.
Then for any $0\le k \le d$, and any convex sets $K, L \subset M_\R$ with non-empty interiors, 
\begin{equation}\label{eq:cvxbound}
\Vol \left( \A ^n(K) [k], L[d-k]\right)    \asymp \|\wedge^k \!\!\A^n\|, 
\end{equation}
where $\wedge^k\!\A^n$  denotes the natural induced linear map on $ \wedge^k M_\R$, and 
$\|\cdot\|$ is any norm on $\End(\wedge^k M_\R)$. 
\end{thm}

It remains to prove Theorem ~\ref{thm:cvx1}. We first present a simple proof in the case when $\A$ is diagonalizable over $\R$ in Section ~\ref{sec:diag}. To deal with the general case, we rely on the Cauchy-Crofton
formula. Some basic material on the geometry of the affine Grassmannian is given in Section  ~\ref{sec:grass}, and the proof of Theorem~\ref{thm:cvx1} is then given in~Section  ~\ref{sec:cvx1}.

Note that, since all norms on $\End(\wedge^k M_\R)$ are equivalent, it
suffices to prove \eqref{eq:cvxbound} for one particular choice of $\|\cdot\|$.

%%%

\subsection{Proof of Theorem~\ref{thm:cvx1} in the diagonalizable case}\label{sec:diag}

Assume that $\A$ is diagonalizable over $\R$, and denote by  
 $\rho_1,\ldots, \rho_d$ its eigenvalues,  ordered so 
that $|\rho_1|\geq \ldots \geq |\rho_d|$.

Let us  first compute $\|\wedge^k \!\A^n\|$.
We fix a basis  $e_1, \ldots ,e_d$  of $M_\R$ that diagonalizes $\A$ so that
$\A e_j = \rho_j e_j$ for all $j$.
For any $k$-tuple $I=\{i_1,\ldots, i_k\}$ of distinct elements in $\{1, \ldots , d\}$, we write 
$e_I := e_{i_1}\wedge \cdots  \wedge e_{i_k}$ and $\rho^I := \prod_1^k \rho_{i_j}$. Then $(\wedge^k \!\A) (e_I) = \rho^I \, e_I$ and 
the collection of $e_I$'s forms a basis of $\wedge^k M_\R$ that diagonalizes $\wedge^k\!\! \A $.
If $\|\cdot\|_{\textrm{sup}}$ is the supremum norm with respect to this basis, then 
\begin{equation*}
\|\wedge^k \!\!\A^n\|_{\textrm{sup}} =  \prod_{j=1}^k |\rho_j|^n. 
\end{equation*}

We now turn to the computation of the mixed volume $\Vol ( \A^n(K) [k], L[d-k])$. First, fix a Euclidean metric $g$ on $M_\R$ such that the basis $e_1,\ldots e_d$ is orthonormal, and let $\Vol_g$ denote the induced volume element. Then there exists a constant $C>0$ such that  $\Vol_g(K)  = C\, \Vol (K)$ for any convex body $K\subset M_\R$.
 It follows that 
$$\Vol \left ( \A^n(K) [k], L[d-k]\right )= C^{-1} \Vol_g \left ( \A^n(K) [k], L[d-k]\right ). $$

Since $K$ and $L$ have non-empty interiors, by arguments as in the proof of Proposition ~\ref{prop:compare}, 
one can show that 
\begin{equation}\label{byta}
\Vol_g \left ( \A^n(K) [k], L[d-k]\right ) \asymp \Vol_g \left ( \A^n(K') [k], K'[d-k]\right ),
\end{equation}
where $K'$ is any convex set with non-empty interior. 

We will compute the right hand side of \eqref{byta} when $K'$ is a polydisk. 
For $r=(r_1,\ldots, r_d)\in\R^d_{\geq 0}$, let $\mathbb D_r$ be the polydisk $\mathbb{D}_r:=\{\sum x_j e_j,\, |x_j|\leq r_j/2\}\subset M_\R$. Note that  $t \mathbb{D}_r+ \tau \mathbb{D}_s=\mathbb{D}_{ t r+\tau s}$ for $r, s \in \mathbb R^d_{\geq 0}$ and $t ,\tau\in \R_{\geq 0}$. It follows that  
$\Vol_g(t \mathbb{D}_r+\tau \mathbb{D}_s)=\prod_{j=1}^d(t  r_j+\tau s_j)$. Thus by \eqref{minkowski}, 
$\Vol_g (\mathbb{D}_r[k], \mathbb{D}_s[d-k]) 
={\binom{d}{k}}^{-1} \sum r^I s^{I^C}$, 
where the sum runs over all multi-indices $\{i_1,\ldots, i_k\}\subset
\{1,\ldots, d\}$, $r^I:= \prod_1^kr_{i_j}$, and $I^C:=\{1,\ldots, d\}\setminus I$. 
Let ${\bf 1}:= (1, \ldots , 1) \in \R_{\geq 0}$. 
Then $\A^n \mathbb D_\1= \mathbb D_{(|\rho_1|^n,\ldots, |\rho_d|^n)}$, and 
\begin{equation*}
\Vol_g\left (\A^n(\mathbb{D}_\1)[k], \mathbb{D}_\1[d-k]\right ) = {\binom{d}{k}}^{-1} \sum_{|I| = k} |\rho^I|^n \asymp \max_I |\rho^I|^n = \prod_{j=1}^k|\rho_j|^n. 
\end{equation*}
This concludes the proof of Theorem ~\ref{thm:cvx1} in the diagonalizable case. 
\begin{rmk}
The same idea can be used to treat the case $A$ is diagonalizable over $\C$. In this case, one chooses $K'$ to be an adequate product of segments and two-dimensional disks.
\end{rmk}

%\bigskip
%
%
%Let us now prove that $\lambda_k(\phi_A) = \prod_1^k |\rho_j|$. We already know it is true when $A$ is diagonalizable over $\C$. We may now pick a sequence of matrices with real coefficients that are all diagonalizable over $\C$, and such that  $A_l \to A$. For each $l$, denote by $ |\rho_1^l| \ge ... \ge |\rho^l_d|$ the eigenvalues of $A_l$, so that  $\lambda_k(\phi_{A_l}) = \prod_1^k |\rho_j^l|$.  
%Since $\wedge^k A_l \to \wedge^k A$, by taking the spectral norm, we get $\prod_1^k |\rho_j^l| \to \prod_1^k |\rho_j|$. To conclude, we need to show that $\lambda_k(\phi_{A_l})\to \lambda_k(\phi_A)$.
%
%Note that the commutator $[A, A_l] \to 0$. For any $\e>0$, we thus have 
%$ \| (A - A_l)(v) \| \le \e\, \| v \|$, and  $\| [A, A_l] \| \le \e\, \| v\|$ for all $v$ and all $l$ large enough.
%We fix $\rho$, such that $\sup \|A_l \| \le \rho$.
%Let us prove by induction that 
%$$
%A^n(\B_1) \subset A_l^n (\B_1) + \B_{\e \rho_n}
%$$
%with $\rho_{n+1} =(n+1) \e \rho^n + \e\rho \rho_n$.
%We first have  the inclusion $A(\B_1) \subset A_l (\B_1) + \B_{\e}$, so that we can take
%$\rho_1 = \rho$. 
%$$
%A^{n+1}(\B_1)
%\subset
%A( A_l^n (\B_1) + \B_{\e\rho_n})
%\subset 
%A_l^n A (\B_1) + \B_{n\e\rho^{n}} + \B_{\e \rho\rho_n}
%\subset 
%A_l^{n+1} (\B_1) + \B_{(n+1)\e\rho^n} + \B_{\e \rho\rho_n}
%$$
%

%%%%

\subsection{The affine Grassmannian}\label{sec:grass}
For $k=1, \ldots, d-1$, we denote by $\gr(d,k)$ the Grassmannian of linear subspaces of $M_\R$ of dimension $k$, and by $\graff(d,k)$ the Grassmanian of affine $k$-dimensional subspaces. Then $\gr(d,k)$ and $\graff(d,k)$ are smooth manifolds, and 
there is a natural projection map $\varpi : \graff(d,k) \to \gr(d,k)$ sending
an affine subspace to the unique linear subspace that is parallel to it. 
The preimage $\varpi^{-1} (H)$ of $H\in\gr(d,k)$ is canonically
identified with $M_\R/H$, and hence we can view $\graff(d,k)$ as the total space of a rank $d-k$ vector bundle over $\gr(d,k)$. 
For $v\in M_\R$, and $H\in\gr(d,k)$, we write $v+H\in\graff(d,k)$ for the affine space obtained by translating $H$ by $v$, so that  
$\varpi(v+H) = H$. Note that the zero section of $\varpi$ is the natural inclusion map $\gr(d,k)\hookrightarrow \graff(d,k)$ given by viewing a linear space as an affine one.

The tangent space $T_H\gr(d,k)$ is canonically isomorphic to $\Hom (H,M_\R/H)$, see~\cite[Ex.~VI.4.1.3]{Sh}. 
It follows that at any point $v+H$ in the affine Grassmannian, we have the exact sequence
\begin{equation}\label{sequence}
0 \to M_\R/H \to  T_{v+H} \graff(d,k) \to \Hom (H,M_\R/H) \to 0~.
\end{equation}

Suppose we are given an invertible \emph{affine map} $\A_\mathrm{aff}: M_\R \to M_\R$, with 
linear part $\A$, i.e.,  $A_{\aff}=A+w$ for some $w\in M_\R$. 
Then $A_\mathrm{aff}$ induces 
smooth maps $A_\mathrm{aff}:
\graff(d,k)\to\graff(d,k)$ and $A: \gr(d,k)\to\gr(d,k)$. 
For any tangent vector $\tau\in T_H\gr(d,k)$, interpreted as a linear
map $\tau : H\to M_\R/H$, we have 
$$ dA_H (\tau) =  A \circ \tau \circ A^{-1}\in T_{A(H)}\gr(d,k).$$
The differential of $A_\mathrm{aff}$ at $v+H\in\graff(d,k)$ is computed analogously using \eqref{sequence}.

\smallskip

Let us fix a Euclidean metric $g$ on $M_\R$, and denote by $\Vol_g$ the
induced volume element. Note that $\varpi^{-1}(H)\simeq M_\R/H$ is canonically identified with $H^\perp$. 
We will see that there are natural induced Riemannian metrics on $\gr(d,k)$
and $\graff(d,k)$. 
First, note that there is a natural action of  the orthogonal group $\rO(M_\R) \simeq \rO(d)$ on $\gr(d,k)$ 
sending $(\phi,H)$ to $\phi(H)$. Since the stabilizer of  a $k$-dimensional
subspace $H\subset M_\R$ is 
$\rO(H)\times\rO(H^\perp) \simeq \rO(k)\times \rO(d-k)$, the
Grassmanian $\gr(d,k)$ 
is diffeomorphic to the homogeneous space $\rO(d)/(\rO(k)\times
\rO(d-k))$. 
There is a Riemannian metric on  $\rO(d)$ that is both left and right invariant by the action of $O(d)$; it is given  
by the pairing $(X,Y) \mapsto - \Tr (XY)$ in its Lie algebra. This
metric induces a  Riemannian metric $g_{\gr}$ on
$\gr(d,k)$ invariant by the action of $\rO(d)$,
see~\cite[Thm.~2.42]{GHL}.

We saw above that $\varpi: \graff(d,k)\to \gr(d,k)$ identifies
$\graff(d,k)$ as the total space of a vector bundle over $\gr(d,k)$.
Any fixed  affine $k$-plane  has a canonical representation $v+ H$ with $v \in H^\perp$. The section $\gr(d,k) \to\graff(d,k)$ sending $H'$ to  $v + H'$ determines 
a subspace $\xi_{v+H}$ in the tangent space of $\graff(d,k)$ at $v+H$ such that 
$d\varpi: \xi_{v+H} \to T_H \gr(d,k)$ is an isomorphism, and
$T_{v+H} \graff(d,k) = \ker d\varpi \oplus\xi_{v+H}$.
We may therefore endow $\graff(d,k)$ with the unique  metric $g_{\graff}$ making  this decomposition orthogonal, such that the restriction $d\varpi: \xi_{v+H}\to T_H \gr(d,k)$
and  the  isomorphism $M_\R /H \simeq H^\perp$ are isometries.
Note that this metric is invariant by $\rO(d)$ but \emph{not} by translations.
However, any translation preserves the fibers of $\varpi$ and their
restriction to each fiber is an isometry.

\smallskip

Recall that any Riemannian metric $h$ on a manifold $M$ defines a volume
element $\Vol_h$ 
on $M$ 
(and a volume form $\dVol_h$ on $M$ if it is oriented),  
see~\cite[Sect.~2.7]{GHL}. If $x=(x_1,\ldots, x_s)$ are local
coordinates on $M$, then locally $h=\sum h_{ij}dx_i\otimes dx_j$ and
$\dVol_h=\sqrt{|\det(h_{ij})|}\, |dx_1 \wedge ... \wedge dx_s|$. 
We will denote by $\Vol_{\gr}$ and $\Vol_{\graff}$ the volume elements
defined by the metrics $g_{\gr}$ and $g_{\graff}$, respectively.  
In fact, $\Vol_{\gr}$
is the unique (up to a scaling factor) volume element that is invariant under the
action of $\rO(d)$.

Recall that an affine map is an \emph{affine orthogonal
  transformation} if (and only if) its linear part is orthogonal. 
\begin{prop}\label{prop:fubini}
The volume element $\Vol_{\graff}$ on $\graff(d,k)$ is invariant by
the group of all affine orthogonal transformations. Moreover it satisfies the 
following Fubini-type property:

\begin{equation}\label{eq:fubini}
\int_{\graff(d,k)} h\, \dVol_{\graff} = \int_{H\in \gr(d,k)} 
\left(
\int_{H^\perp}\, h\, \dVol_{g|_{H^\perp}} \right)\, \dVol_{\gr},
\end{equation}
for any Borel function $h$ on $\graff(d,k)$.
\end{prop}

\begin{proof}
Let us first prove~\eqref{eq:fubini}.  Pick $H\in\gr(d,k)$ and a neighborhood $H\in\mathcal U\subseteq \gr(d,k)$ with local coordinates $y=(y_1,\ldots, y_{k(d-k)})$.  
Then locally in $\mathcal U$, $g_{\gr}$ is of the form $g_{\gr}=\sum b_{ij}(y)\, dy_i\otimes dy_j$. 
Moreover we can choose a local trivialization $\mathcal U\times \R^{d-k}$ of $\graff(d,k)\to\gr(d,k)$ with coordinates $(x,y)$, where $\varpi(x,y) = y$ and $x=(x_1,\ldots, x_{d-k})$ are coordinates in $\R^{(d-k)}\simeq H^\perp$. 
Since  $\varpi:\graff(d,k)\to\gr(d,k)$ is a Riemannian submersion, and the restriction of $g_{\graff}$ to the fiber $\varpi^{-1}(H)\simeq H^\perp$ is the constant metric $g|_{H^\perp}=\sum a_{ij}(H) \, dx_i\otimes dx_j$, locally in the trivialization, $g_{\graff}(x,y) =  \sum a_{ij}(y) \,dx_i\otimes dx_j + 
\sum b_{ij}(y)\, dy_i\otimes dy_j$. 
Consequently 
$d\Vol_{\graff} =  \sqrt{\det(a_{ij}(y))} \,  \sqrt{\det(b_{ij}(y))}\,dx\wedge dy$. After a partition of unity we may assume that $h$ has support in a small neighborhood of $v+H\in\graff(d,k)$ and thus
$$
\int h\, \dVol_{\graff(d,k)}  = 
\int   \left(  \int h \,\sqrt{\det(b_{ij}(y))} \, dx \right) \wedge \sqrt{\det(a_{ij}(y))} \, dy,
$$
which proves \eqref{eq:fubini}. 

To prove the first part of the proposition we need to show that $\Vol_{\graff}$ is invariant under linear orthogonal transformations and translations. First, since $g_{\graff}$ is invariant under $O(d)$, so is $\Vol_{\graff}$. Next, let $A_w:=\id + w: \graff(d,k)\to \graff(d,k)$ be the translation by $w\in M_\R$. Since $g$, and thus $g|_{H^\perp}$, is invariant under translations on $H^\perp$, 
$$
\int_{v\in H^\perp}\, h\circ A_w(v)\, \dVol_{g|_{H^\perp}}= \int_{v\in H^\perp}\, h(v+w')\, \dVol_{g|_{H^\perp}} = \int_{v\in H^\perp}\, h(v)\, \dVol_{g|_{H^\perp}},
$$
where $w'$ is the orthogonal projection of $w$ onto $H^\perp$. Thus, using  
\eqref{eq:fubini}, $\int (h \circ A_w)\,  \dVol_{\graff}= \int h\,  \dVol_{\graff}$; this shows that $\Vol_{\graff}$ is invariant under translations. 
\end{proof}

Let $A:\gr(d,k)\to\gr(d,k)$ and $A_{\mathrm{aff}}: \graff(d,k)\to\graff(d,k)$ be the maps induced by an invertible affine map $A_{\mathrm{aff}}:M_\R\to M_\R$ with linear part $A$. Recall that the \emph{Jacobians} of $A$ and $A_{\aff}$, respectively, are the unique smooth functions 
$J A: \gr(d,k)\to \R_{\geq 0}$ and $J A_{\mathrm{aff}} :\graff(d,k)\to\R_{\geq 0}$
that satisfy the change of variables formula, i.e., 
\begin{eqnarray}
\int_{\gr(d,k)} h \, \dVol_{\gr} 
&=&
 \int_{\gr(d,k)} (h\circ A) \, J A \, \dVol_{\gr},
\label{eq:chg1}\\
\int_{\graff(d,k)} h \, \dVol_{\graff} 
& =&  \int_{\graff(d,k)} (h\circ A_\aff) \, JA_\mathrm{aff} \, \dVol_{\graff}
\label{eq:chg2}
\end{eqnarray}
for any integrable functions $h$ on $\graff(d,k)$ and $\gr(d,k)$, respectively.

Given a linear map $A: M_\R\to M_\R$, let $\Phi_\A : \gr(d,k)\to \R_{\geq 0}$
be the map that maps $H$ to (absolute value of) the Jacobian of the induced linear map $A : M_\R/H \to M_\R/A(H)$, computed with respect to the volume elements $\Vol_{g|_{H^\perp}}$ and $\Vol_{g|_{\A(H)^\perp}}$ defined by $g$ on $M_\R/H\simeq H^\perp$ and $M_\R/A(H)\simeq A(H)^\perp$, respectively. In other words, $\Phi_\A $ is the unique function  that satisfies 
\begin{equation}\label{sweetwaters}
\int_{M_\R/A(H)} h \, \dVol_{g|_{\A(H)^\perp}} =\int_{M_\R/H} (h \circ A)\, \Phi_\A(H) \, \dVol_{g|_{H^\perp}}
\end{equation}

The following is a key lemma in the proof of Theorem ~\ref{thm:cvx1}.

\begin{lem}\label{prop:jacobian-graff}
Let $A:M_\R\to M_\R$ be any invertible linear map. 
Then for any linear space $H\subset M_\R$, and any $v\in M_\R$, we have
\begin{equation*}\label{eq:jac-graff}
J A_\mathrm{aff} (v+H) =
J A (H)\times \Phi_\A (H).
\end{equation*}
\end{lem}

\begin{proof}
Recall that, by the definition of $g_{\graff}$, the tangent space of
$\graff(d,k)$ at $v+H$
splits as an orthogonal sum $M_\R/H\oplus T_H\gr(d,k)$. The differential of $A_\aff$
does not preserve this orthogonal decomposition in general but sends $M_\R/H$ to 
$M_\R/A(H)$; the tangent space at $A_\aff(v+H)=A_\aff(v)+A(H)$
orthogonally splits as $M_\R/A(H)\oplus T_{\A(H)}\gr(d,k)$. 
Choose (local) orthonormal bases $v_j$, $w_j$, $v_j'$, and $w_j'$ of $M_\R/H\simeq H^\perp$,
$T_H\gr(d,k)$ (at $H$),
$M_\R/A(H)\simeq A(H)^\perp$ and
$T_{\A(H)}\gr(d,k)$ (at $A(H)$), respectively. 
Then, in light of \eqref{eq:chg2}, the Jacobian of $\A$ at $v+H$ is the
absolute value of the determinant of the matrix $dA$ with respect to the bases $\{v_j,
v_j'\}$ at $v+H$ and $\{w_j, w_j'\}$ at $A(v) + A(H)$. This matrix,
however, is block diagonal and so its determinant is the product of two
determinants: one of the matrix $dA:M_\R/H\to M_\R/A(H)$ with respect
to the bases $\{v_j\}$ and $\{w_j\}$ and one of the matrix of
$dA:T_H\gr(d,k)\to T_{A(H)}\gr(d,k)$ in the bases $\{v_j'\}$ and
$\{w_j'\}$. In light of \eqref{sweetwaters} and \eqref{eq:chg1}, this concludes the proof. 
\end{proof}

%%%

\subsection{Proof of Theorem~\ref{thm:cvx1}}\label{sec:cvx1}
Let $\bB\subset M_\R$ be the unit ball
with respect to the metric $g$ on $M_\R$ from last section. Then, by arguments as in Section ~\ref{sec:diag},
\begin{equation}\label{snyta}
\Vol\left ( \A^n(K) [k], L[d-k]\right)\asymp \Vol_g \left ( \A^n(\bB) [k], \bB[d-k]\right ).
\end{equation}
To compute the right hand side of \eqref{snyta} we will apply 
the Cauchy-Crofton formula, see~\cite[formula 4.5.10]{S}, which asserts that there exists a universal
constant $C>0$ such that  for any convex set $K\subset M_\R$:
\begin{equation}\label{eq:cauchy-crofton}
\Vol_g \left( K [k], \bB[d-k]\right) = C\,\Vol_{\graff} \left\{
v + H\in \graff(d,d-k),\, (v+ H)\cap  K \neq \emptyset
\right\}~,
\end{equation}
where $\Vol_{\graff}$ is defined as in Section ~\ref{sec:grass}.

Now 
\begin{multline*}
\frac1C \Vol_{\graff} \left( \A^n(\bB) [k], \bB[d-k]\right) 
= 
\int_{v+H\in\graff(d,d-k), \, (v+H) \cap \A^n(\bB)\neq\emptyset} \dVol_{\graff}
=\\
\int_{v+H\in\graff(d,d-k), \, (v+H) \cap \bB\neq \emptyset} \, J \A^n (v+H)\, \dVol_{\graff} 
=\\
\int_{v+H\in\graff(d,d-k), \, (v+H) \cap \bB\neq \emptyset}  (J \A^n \times \Phi_{\A^n})(H)\, \dVol_{\graff}
=\\
\int_{H\in\gr(d,d-k)} \left (\int_{v\in H^\perp, (v+H)\cap \bB\neq \emptyset} \dVol_{g|_{H^{\perp}}}\right )\, 
(J \A^n \times \Phi_{\A^n})(H)\, \dVol_{\gr}
=\\
V_{k} \int_{H\in\gr(d,d-k)} \Phi_{\A^n}  \times J \A^n \, \dVol_{\gr} 
= 
V_{k} \int_{\gr(d,d-k)} (\Phi_{\A^n}\circ \A^{-n}) \, \dVol_{\gr},
\end{multline*}
where $V_{k}$ is the volume of the orthogonal projection of $\bB$ onto $H^\perp$, i.e., the volume of the standard $k$-dimensional ball in Euclidean space, and $\Vol_{\gr}$ and $\Phi_{\A}$ are defined as in Section ~\ref{sec:grass}. Here we have used \eqref{eq:cauchy-crofton}, \eqref{eq:chg2}, Lemma ~\ref{prop:jacobian-graff}, \eqref{eq:fubini}, and \eqref{eq:chg1} for the first, second, third, fourth, and last equality, respectively.

To sum up, 
\begin{equation}\label{sumup}
\Vol \left( \A^n(K) [k], L[d-k]\right) 
\asymp 
\int_{\gr(d,d-k)}
(\Phi_{\A^n}\circ \A^{-n}) \,
\dVol_{\gr}, 
\end{equation} 
We will prove Theorem ~\ref{thm:cvx1} by estimating the right hand side of \eqref{sumup}. For that we will need the following two lemmas.

\begin{lem}\label{lem:jaco}
Let $H\subset M_\R$ be a linear subspace of codimension $k$ and let
$A: M_\R\to M_\R$ be a linear map with $\det (\A)\neq 0$. 
Then for any $\gamma\in\wedge^{d-k} M_\R$ defining $H$ in the sense
that $\gamma \wedge v =0$ if and only if
$v\in H$,
we have
\begin{equation*}
\Phi_{\A} \circ \A^{-1} (H) = |\det (\A)|\, \frac{\|\wedge^{d-k}\A^{-1} (\gamma)\|}{\|\gamma\|}~,
\end{equation*}
where $\wedge^{d-k} A^{-1}$ is the induced linear map on
$\wedge^{d-k}$ and  $\|\cdot\|$ is the natural norm on  $\End(\wedge^{d-k} M_\R)$ induced by $g$.
\end{lem}

\begin{lem}\label{lem:lower}
Let $(V, \|\cdot\|)$ be a finite dimensional normed vector space, and let 
$h: V \to V$ be a linear map with $\det (h)\neq 0$. 
Moreover, for $v\in V$, let 
$$
\tau_h(v) := \inf_n \frac{\| h^n( v)\|}{\|v \|\, \|h^n\|}. 
$$
Then $\tau_h : V \to \R_{\geq 0}$ is an upper semicontinuous function and $\{ \tau_h = 0\}$ is a 
proper linear subspace of $V$. 
\end{lem}

Set $ h := \det(A)\, (\wedge ^{d-k} A^{-1})$. Observe that, since the pairing $\wedge ^k M_\R \times
\wedge^{d-k} M_\R \to \wedge^d M_\R$ is perfect, in fact, $\|h^n\|= \|\wedge^k A^n\|$.

For any $H \in \gr(d,d-k)$, we pick $\gamma(H)\in\wedge^{d-k}M_\R$ defining it. 
Then $\gamma(H)$ is unique up to a scalar factor, and the induced map  $\pl: \gr(d,d-k) \to \P(\wedge^{d-k}M_\R)$ is the Pl\"ucker embedding of  $\gr(d,d-k)$. 
Lemma ~\ref{lem:jaco} can be rephrased as follows:
\begin{equation}\label{slutet}
\Phi_{A^n}\circ A^{-n} (H) = \frac{|h^n(\pl(H))|}{|\pl(H)|}.
\end{equation}
Note that the right hand side of \eqref{slutet} is well defined by homogeneity.
The image of $\gr(d,d-k)$  under $\pl$ is not contained in any proper linear
subspace of $\wedge^{d-k}M_\R$. Therefore, by Lemma ~\ref{lem:lower}, there
is a non-empty open set  $\mathcal U\subset \gr(d,d-k)$, such that
$\tau_h $ restricted to $\mathcal U$ is strictly positive. 
In particular, $\mu := 
\int_{\mathcal U} \tau_h(\pl(H)) >0$. 
Consequently,  
\begin{equation}\label{nedre}
\int_{\gr(d,d-k)} (\Phi_{\A^n}\circ \A^{-n}) \, \,  \dVol_{\gr} 
\geq 
\| h^n\| \int_{\gr(d,d-k)} \tau_h (\pl(H)) \dVol_{\gr} 
\geq 
 \mu \, \|\wedge^k A^n\|~.
\end{equation}
Now \eqref{eq:cvxbound} follows from \eqref{sumup},  \eqref{nedre}, and the trivial upper bound
$\Phi_{A^n} \circ A^{-n} \le \| h^n\|=\|\wedge^k A^n\|$.
Thus we have proved Theorem ~\ref{thm:cvx1}.

It remains to prove the lemmas. 
\begin{proof}[Proof of Lemma~\ref{lem:jaco}] 
Pick $H\in \gr(d,d-k)$. Choose orthonormal bases $e_1,\ldots, e_d$ and $f_1,\ldots, f_d$ of $M_\R$ such that 
$(e_1, \ldots , e_k)\in \A^{-1}(H)^\perp\simeq M/\A^{-1}(H)$ and  $(f_1, \ldots , f_k)\in H^\perp\simeq M/H$. Then $\A=\sum a_{ij} f_i\otimes e_j^*$ for some $a_{ij}\in \R$, and $\Phi_\A\circ \A^{-1}(H)$ is by definition equal to 
$\left| \det (a_{ij})_{1\le i,j \le k}\right|$. 
On the other hand the vector $\gamma = f_{k+1} \wedge \cdots \wedge f_d$ defines $H$,
and $$\wedge^{d-k}\A^{-1} (\gamma) =  \frac{e_{k+1}\wedge \cdots\wedge e_d}{\det (a_{ij})_{k+1\le i,j \le d}} 
 =\pm \frac{\Phi_\A\circ \A^{-1}(H)}{|\det(A)|}\, e_{k+1}\wedge \cdots\wedge e_d 
 ~.$$ 
We conclude noting that $|e_{k+1}\wedge \cdots \wedge e_d|=|\gamma|=1$.
\end{proof}

\begin{proof}[Proof of Lemma~\ref{lem:lower}]
For each $n$ the function $v \mapsto  \frac{\| h^n( v)\|}{\|v \|\, \|h^n\|}$ is continuous, and so $\tau_h$ is the infimum of a sequence of continuous functions, which implies that it is upper semicontinuous.

Let us now describe the zero locus of $\tau_h$. First, assume that
there are no non-trivial subspaces of $V$ that are invariant under
$h$. Choose a basis of $V$
such that the matrix (also denoted by $h$) of $h$ is in Jordan normal
form (over $\C$), and let $x_1,\ldots, x_d$ with $d = {\dim V}$ be the corresponding
coordinates. Then 
\begin{equation}\label{matrisrepr}
h =
\left[
\begin{array}{ccccc}
\rho & 1 & 0 & \ldots & 0
\\
0  & \rho  & 1 & \ldots & 0
\\
0 & 0 & \rho & \ldots & 0
\\
&& \ddots &&
\\
0 & 0 & 0 & \ldots & 1 
\\
0 & 0 & 0 & \ldots & \rho
\end{array}
\right] 
\text{ and } 
h^n \asymp
\left[
\begin{array}{ccccc}
\rho^n  & n\rho^n & n^2\rho^n & \ldots & n^{d-1}\rho^n
\\
0  & \rho^n  &  n\rho^n & \ldots & n^{d-2}\rho^n
\\
0 & 0 & \rho^n & \ldots & n^{d-3}\rho^n
\\
&& \ddots &&
\\
0 & 0 & 0 & \ldots & n \rho^n
\\
0 & 0 & 0 & \ldots & \rho^n
\end{array}
\right].
\end{equation}
In this case,  $\{\tau_h =0\}$ is precisely the hyperplane $\{x_{\dim\!
  V}=0\}$. 

In the general case, we decompose $V=\bigoplus  W_i$ into minimal
$h$-invariant subspaces $W_i$. 
Let $\rho_i$ be the modulus of the unique eigenvalue
of $h|_{W_i}$ and write $d_i := \dim (W_i)$. 
Set $\rho := \max\{\rho_i\}$, 
$I_0 := \{i,\, \rho_i = \rho\}$, $\d := \max\{d_i, \, i\in I_0\}$, and $I_*:= \{ i\in I_0, d_i = \d\}$.
Then $\tau_h = \max \{ \tau_{h|_{W_i}}\circ p_{W_i}, \, i\in I_*\}$, 
where $p_i:V\to W_i$ is the natural projection, 
and 
$\{\tau_h =0\}$ is the direct sum of the $W_i$ with $i\notin I_*$
and the hyperplanes 
$ \{\tau_{h|_{W_i}}\ =0\}\subset W_i$ for $i\in I_*$; in particular, 
$\{\tau_h=0\}$ is linear, and 
since $I_*$ is non-empty, it is a proper subspace of $V$. 
\end{proof}

\subsection{Proof of Corollary B}\label{pfB}
Choose a basis of $M_\R$ such that the matrix of $A$ is in Jordan
normal form, and let $\|\cdot\|_{\sup}$ be the supremum norm with
respect to the induced basis on $\wedge^k M_\R$. Then $\|\wedge^k A^n
\|_{\sup}\asymp n^r|\rho_1|^n\cdots |\rho_k|^n$ for some integer
$0\leq r\leq d-1$, cf. the proof of Lemma ~\ref{lem:lower}. Hence
$\l_k=\lim_n(\deg_k(\f_A^n))^{1/n}=|\rho_1|\cdots |\rho_k|$.

%%%%%%%%%%

\section{Proof of Theorem~D}\label{sec:pfD}

%%%%%

As well as Theorem~A, Theorem~D can be proved by controlling the
growth of mixed volumes under the linear map $A:M_\R\to M_\R$.

We fix a euclidean metric $g$ on $M_\R$.
Given a subspace $H\subset M_\R$, let $\Vol_H$ denote the volume element 
on $H$ induced by $g$. 
Moreover, let $p_{H}$ denote the orthogonal projection onto $H$.

\begin{thm}\label{thm:fine-mixed}
Let $\A : M \to M$  be  a group morphism such that $\det(\A) \neq 0$, 
with eigenvalues
$|\rho_1| \ge \ldots  \ge |\rho_d|$.  
Suppose that $\kappa := |\rho_{k+1}| / |\rho_{k}|<1$, and write
$V_u:= \oplus_{i\le k} \ker (\A-\rho_i \id)^d$, and $V_s:= \oplus_{i> k} \ker (\A-\rho_i \id)^d$. Then there exists an integer $r\ge0$, such that 
for any two (non-empty) convex sets $K, L\subset M_\R$, 
\begin{equation}\label{eq:current1}
 \frac1{\l_k^n} \Vol \left(\A^n (K) [k], L[d-k]\right)  
 = 
 \Vol_{V_u} \left(p_{V_u/V_s}(K)\right) \,  \Vol_{V_u^\perp} \left (
   p_{V_u^\perp}(L) \right )
+ \O (n^r\, \kappa^n) ~,
\end{equation}
where $p_{V_u/V_s}$ denotes the projection onto $V_u$ parallel to
$V_s$.  
\end{thm}

Note that, by Corollary B, the condition $\kappa<1$ is equivalent to
\eqref{strikt}. Recall from Proposition~\ref{volume} that $\deg_k(\f_\A^n)=
d!\Vol (\A^n (\Sigma_d)[k], \Sigma_d [d-k])$, where $\Sigma_d$ is the
standard simplex. Thus, noting that $\kappa\l_k=(\l_{k+1}\l_{k-1})/\l_k$, 
Theorem ~\ref{thm:fine-mixed} gives \eqref{lor} with 
$C= d! \Vol (p_{V_u/V_s}(\Sigma_d)) \,  \Vol
(p_{V_u^\perp}(\Sigma_d)) >0$. Taking Theorem ~\ref{thm:fine-mixed} for granted 
this concludes the proof of Theorem~D.

\begin{rmk}\label{general}
Note that Theorem \ref{thm:fine-mixed} applied to $K=L=P_D$, under the
assumption in Theorem D, gives the following version of \eqref{lor}:
$$\deg_{D,k}(\f_A^n)  = C \l_k^n + \O \left(n^r\, \left(\frac{\l_{k-1} \, \l_{k+1}}{\l_k}\right)^{n}\right),$$
where $C= d! \Vol (p_{V_u/V_s}(P_D)) \,  \Vol
(p_{V_u^\perp}(P_D))> 0$. 
\end{rmk}

%%%%%

\subsection{Proof of Theorem~\ref{thm:fine-mixed}}

For the proof  
we will need the
following two lemmas on mixed volumes.

\begin{lem}\label{lem:flat}
Let $H \subset M_\R$ be subspace of dimension $k$. 
Then for any convex sets $L_1, \ldots, L_{d-k},$  and $K$ in $M_\R$
such that $K \subset H$, 
\begin{equation}\label{eq:flatvol}
\Vol_{M_\R} \left (K[k], L_1, \ldots , L_{d-k}\right ) = \Vol_H (K) \, \Vol_{H^\perp}\left (p_{H^\perp}(L_1), \ldots, p_{H^\perp}(L_{d-k})\right ).
\end{equation}
\end{lem}

\begin{proof}[Proof of Lemma~\ref{lem:flat}]
Recall that we have the following polarization formula:
$$(d-k)! \Vol(K[k], L_1, \ldots , L_{d-k}) = \sum_{I\subset \{ 1, \ldots , d-k\}} (-1)^{(d-k- |I|)} \Vol (K[k], (\sum_I L_i)[d-k])~.$$
By multilinearity, we may thus assume that $L_1 = \ldots  = L_{d-k}=L$. Fix $t\in\R$. Then,
by Fubini's theorem, 
\begin{equation}\label{snart}
\Vol ( tK  + L) = \int_{v\in H^\perp} \Vol_H \left ( (tK + L) \cap
  (v + H) \right ) \, \dVol_{H^\perp},
\end{equation}
where we have identified the volume element induced by $\Vol$ on $v+H$  
with $\Vol_H$. 
Let $L_v := L \cap (v + H)$. 
Since $K$ is included in $H$, we have 
$(tK + L) \cap (v + H)  = tK + L_v $, and so the right hand side of
\eqref{snart} equals $\int_{v\in H^\perp} \Vol_{H} ( tK + L_v  ) \,
\dVol_{H^\perp}$.
Note that $L_v \neq \emptyset$ if and only if $v \in p_{H^\perp}(L)$, in which case $\Vol_H ( tK + L_v)=  t^k \Vol(K) + \O(t^{k-1})$. We conclude that 
$$
\Vol (tK  + L) = \int_{p_{H^\perp}(L)} \left(t^k \Vol(K) +\O(t^{k-1}) \right)\,
\dVol_{H^\perp} 
= t^k \Vol(K) \Vol_{H^\perp}(L)
+ \O(t^{k-1}),
$$ 
which implies~\eqref{eq:flatvol}.
\end{proof}

Recall that the \emph{Hausdorff distance} $d_H(K,L)$ between two (non-empty) sets $K, L\subset M_\R$ is the infimum of all $\epsilon\in\R_{\geq 0}$ such
that $K\subset L+\bB_{\epsilon}$ and $L\subset K+\bB_{\epsilon}$, where
$\bB_r\subset M_\R$ is a ball with radius $r$ and center $0$.
In the sequel, we write $\bB = \bB_1$.

\begin{lem}\label{lem:continuity}
Let $L_1, \ldots, L_{d-k}$ be convex sets in $M_\R$.
Then there exists a constant $C>0$ such that for any (non-empty) convex sets
$K, K'\subset M_\R$, one has 
\begin{multline}\label{gul}
|\Vol \left (K[k], L_1, \ldots , L_{d-k}\right ) - \Vol \left (K'[k], L_1, \ldots , L_{d-k}\right )|
\le \\ 
C 
\max_{j=1}^{k} \{d_H(K,K')^j \, \Vol \left (K[k-j], \bB[d-k+j]\right ) \}.
\end{multline}
\end{lem}

\begin{proof}[Proof of Lemma~\ref{lem:continuity}]

To simplify notation, write 
$\Vol \left (\cdots, L_1, \ldots , L_{d-k}\right ) =: \Vol \left
  (\cdots, L_i\right )$,  and $\delta := d_H(K,K')$ so that
$K\subset K'+\bB_\delta $, and $K'\subset K+\bB_\delta$.

Assume first that $\Vol \left ( K[k], L_i\right ) \geq \Vol \left ( K'[k],
  L_i\right )$. 
Using the multilinearity of the mixed volume and \eqref{monotone} we
get: 
\begin{multline*}
\Vol \left ( K[k], L_i\right ) - \Vol \left ( K'[k], L_i\right )
\le
\Vol \left ( (K'+\bB_\delta) [k], L_i\right ) - \Vol \left ( K'[k], L_i\right )
= \\
\sum_{\ell=1}^k {k\choose \ell}\, \delta^\ell\, \Vol \left ( K'[k-\ell],  \bB [\ell], L_i\right )
\le \sum_{\ell=1}^k {k\choose \ell}\, \delta^\ell\, \Vol \left ( (K+\bB_\delta) [k-\ell],  \bB [\ell], L_i\right )
=\\
\sum_{j=1}^k C_j \, \delta^j\, \Vol \left ( K[k-j],  \bB
[j], L_i\right )
\leq 
C \max_{j=1}^k \{ (\delta)^j\, \Vol \left ( K[k-j],  \bB [d-k+j] \right )\}
\end{multline*}
for some constants $C_j,C>0$; for the last inequality we have used
that each $L_i$ is contained in $\bB_{\rho_i}$ for $\rho_i\in\R_{\geq 0}$ large
enough. This proves \eqref{gul} in this case. 

If $\Vol \left ( K'[k], L_i\right ) \geq \Vol \left ( K[k],
  L_i\right )$, then \eqref{gul} follows as above noting that $\Vol \left ( K[k-j],
  \bB [d-k+j] \right ) \asymp \Vol \left ( K'[k-j],  \bB [d-k+j]
\right )$. 
\end{proof}

We are now ready to prove Theorem ~\ref{thm:fine-mixed}. Write $p:=
p_{V_u/V_s}$ to simplify notation.  

First, since $p(K)$, as well as $\A^n\circ p
(K)=p\circ
\A^n(K)$, 
is included in the $k$-dimensional subspace $V_u\subset
M_\R$, 
by Lemma ~\ref{lem:flat}, 
\begin{multline}\label{carpet}
\l_k^n \Vol_{V_u} \left (p(K)\right ) \, \Vol_{V_u^\perp} \left
  (p_{V_u^\perp}(L)\right )
=\\
\Vol_{V_u} \left ( p\o \A^n(K)\right)\, 
\Vol_{V_u^\perp} \left (p_{V_u^\perp}(L)\right) = 
\Vol \left ( p\o \A^n(K)[k], L[d-k]\right ).
\end{multline}
Next, note that there exists a constant $C>0$ and an integer $r\geq 0$ such
that $|\A^n(v)| \le C\, n^r\, |\rho_{k+1}|^n\, |v|$ for all $v\in V_s$;
for example, this can be seen using \eqref{matrisrepr}. In particular, 
$|\A^n(v) - p \o \A^n (v)| \leq  C\, n^r |\rho_{k+1}|^n\, |v|$ for all
$v\in M_\R$, from which 
we infer that
\begin{equation}\label{matta}
d_H(\A^n(K), p\o \A^n(K) ) \le C'\, n^r\, |\rho_{k+1}|^n,
\end{equation}
where $C'=C\max_{v\in K}|v|$. Now, 
applying Lemma ~\ref{lem:continuity} to $K = \A^n(K)$, $K'=p\circ
\A^n(K)$, and $L_i
=L$ for all $i$, and using \eqref{carpet} and
\eqref{matta}, we get 
\begin{multline}\label{laban}
\Vol\left (\A^n(K)[k], L[d-k] \right ) - \l_k^n \Vol_{V_u} \left (p(K)\right ) \, \Vol_{V_u^\perp} \left
  (p_{V_u^\perp}(L)\right )\le \\
C'' \max_{j=1}^{k} \{ n^{jr}\, |\rho_{k+1}|^{jn} \, \Vol \left (\A^n(K) [k-j], \bB[d-k+j]\right)\}
\end{multline}
for some constant $C''> 0$. 
Furthemore, Theorem ~\ref{thm:cvx1} implies that  
$$\Vol \left (\A^n(K) [k-j], \bB[d-k+j]\right )\le C_j n^{r_j} \prod_{i=1}^{k-j}|\rho_i|^n$$
for suitable constants $C_j>0$ and integers $r_j\ge 0$, cf. Section
~\ref{pfB}. 
Since $ |\rho_1|\cdots|\rho_{k-j}| |\rho_{k+1}|^j \le \lambda_k \kappa$
for $1\leq j\leq k$ by Corollary B, 
the right hand side of \eqref{laban} is bounded from above by $C'''
n^{r'}( \lambda_k \kappa)^n$
for some constant $C'''>0$ and some integer $r'\geq 0$, which proves
\eqref{eq:current1}.

%%%%%%%%%%

\section{Invariant classes}\label{sec:inv}
%\subsection{Complements: invariant classes}\label{sec:inv}

In fact, Theorem~\ref{thm:fine-mixed} gives more information 
than Theorem D. Keeping the notation from the beginning of Section ~\ref{sec:pfD}, 
consider the currents $T^-:= [V_u, p_{V_u/V_s}]$ and 
$T^+ := [V_u^\perp, p_{V_u^\perp}]$
of degree $(d-k)$ and $k$, respectively, as defined in Example ~\ref{proj-val}. 
Then for polytopes $P,Q\subset M_\Q$, \eqref{eq:current1} reads
\begin{equation*}
\frac1{\l_k^n} \Vol \left(\A^{n}_*[P] [k], [Q][d-k]\right)  
 = 
\langle T^-,[P]\rangle \,  \langle T^+, [Q] \rangle
+ \O (n^r\, \kappa^n) ,
\end{equation*}
and, by multilinearity, using
\eqref{madison}, we get: 
\begin{thm}\label{thm:invariant}
Let $A$, $\kappa$, $V_s$, $V_u$ be as in Theorem ~\ref{thm:fine-mixed}. 
Then there exists an integer $r$ such that, 
for any $\a \in \Pi_k$, and $\b\in\Pi_{d-k}$, 
\begin{equation*}
\frac1{\l_k^n} 
\Vol ( \A^n_{*} \a\cdot  \b) = 
 {d \choose k}\, \langle T^-, \a\rangle \,  \langle T^+, \b \rangle
+ \O (n^r\, \kappa^n) ~.
\end{equation*}
\end{thm} 
In particular, in the space of currents on $\Pi$, 
the convergence 
$$\frac{1}{\l_k^{n}} \, \A^n_{*} \a \to c(\a)\, T^+$$ 
holds for any class $\a\in\Pi_k$, with $c(\a) = {d\choose k} \langle T^-,
\a\rangle$; here we identify $\gamma\in\Pi$ with $T_\gamma\in\mathcal
C$, cf. Example ~\ref{poly-embedding}. 
By Lemmas~\ref{lem:zero} and~\ref{lem:dual}, 
 $\Vol ( \A^{n}_* \a\cdot  \b) = \Vol (\a\cdot   \A^{n*} \b) $
for $\a\in\Pi_k$ and $\b\in\Pi_{d-k}$,
so that, by duality 
$$\frac1{\l_k^{n}} \, \A^{n*} \b \to c'(\b) T^-$$ for any class $\b\in\Pi_{d-k}$, with 
$c'(\b) = {d\choose k}\, \langle T^+, \b\rangle$.
By Theorem~\ref{thm:equiv}, the currents $T^+$ and $T^-$ induce
classes in the universal cohomology of toric varieties,
$\theta^+\in\hproj^k$ and $\theta^-\in \hproj^{d-k}$, respectively.

\begin{cor}\label{klasser}
Let $A$, $\kappa$, $V_s$, $V_u$ be as in Theorem ~\ref{thm:fine-mixed}. 
Then there exists an integer $r$ such that, 
for any complete simplicial fan $\D$, and any classes
$\omega\in H^{2k}(X(\D))$, $\eta\in H^{2(d-k)}(X(\D))$, 
\begin{equation*}
\frac1{\l_k^n} 
(\f_\A^{n})^* \omega\cdot  \eta = 
{d\choose k}\, (\theta_\Delta^-\cdot \omega) \, (\theta_\Delta^+\cdot \eta)
+ \O (n^r\, \kappa^n) .
\end{equation*}
Moreover, if $L$ is an ample class in some projective toric variety, 
and $\omega = L^k$, respectively $\eta = L^{d-k}$, then 
$\langle \theta^+, \omega\rangle> 0$, respectively $\langle \theta^-, \eta \rangle> 0$.
\end{cor}
In particular, $\frac1{\l_k^n} 
(\f_\A^{n})^* \omega$, regarded as a class in $\hinj^k$, converges towards ${d\choose k} (\theta^-\cdot \omega) \, \theta^+$ and by duality $\frac1{\l_k^n} 
(\f_{\A}^{n})_* \eta$, regarded as a class in $\hinj^{d-k}$ converges towards ${d\choose k} (\theta^+\cdot \eta) \, \theta^-$. 
This result is the analog of~\cite[Corollary~3.6]{BFJ} in the context of  monomial maps but in arbitrary dimensions.

%%%%%%%%%%%%%%%%%%%%%%%%%%%%%%%%%%%%%%%%%%%%%%%%

\def\listing#1#2#3{{\sc #1}:\ {\it #2},\ #3.}

\end{document}